\documentclass[twoside,10pt,a4paper]{amsart}

\usepackage{amsmath}
\usepackage{verbatim} 
\usepackage{epsfig}
\numberwithin{equation}{section}

\begin{document}

\title[Cumulants of the Power Law Logistic Model] 
       {Approximations of Cumulants of the Stochastic Power Law Logistic Model}

\date{\today}
\author{Ingemar N{\aa}sell}
\address{Department of Mathematics \\
       The Royal Institute of Technology \\
      S-100 44 Stockholm, Sweden}
 

\begin{abstract}
Asymptotic approximations of the first three cumulants of the quasi-stationary 
distribution of the stochastic power law logistic model are derived. 
The results are based on a system of ODEs for the first three cumulants. 
We deviate from the classical moment closure approach by determining 
approximations without closing the system of equations. 
The approximations are explicit in the model's parameters, conditions for 
validity of the approximations are given, magnitudes of approximation errors 
are given, and spurious solutions are easily detected and eliminated. 
In these ways, we provide improvements on previous results for this model. 
\end{abstract}

\maketitle


\section{Introduction}

We study a stochastic version of the power law logistic model. 
Logistic models are used as models for the size of a population with density 
dependent growth. 
This means that the net growth rate per individual is a decreasing function 
of the population size. 
The simplest decreasing function is the linear one.
It leads to the classical logistic model, whose deterministic version was studied by Verhulst 
already in 1838. 
The stochastic version of this model has been studied extensively, as shown 
by N{\aa}sell (2017). 
More general deterministic power law logistic models have been discussed by 
Banks (1994) and by Tsoularis and Wallace (2002), while stochastic versions of 
such models are treated by Matis, Kiffe, and Parthasarathy (1998), by Renshaw 
(2011), and by Bhowmick, Bandyopadhyay, Rana, and Bhattacharya (2016).

The model that we deal with is a Markov Chain with finite state space and 
continuous time, and with an absorbing state at the origin. 
It can also be described as a finite-state birth-death process.  
Two model formulations are given in Section 2. 
They are mathematically equivalent, but with different parameters. 
The maximum population size $N$ is an important parameter in the second one of 
the two formulations. 
This formulation is basic for our work, since $N$ serves the role of being the large 
parameter for which asymptotic approximations of various quantities can be derived.

The state variable $X(t)$ of either of the two formulations is interpreted 
as the number of individuals in a population. 
Ultimate extinction of this variable is a fact for this model. 
This means in particular that the stationary distribution of $X(t)$ is 
degenerate with probability one at the origin. 
A related and important random variable $X^Q(t)$ is defined by conditoning 
$X(t)$ on non-extinction. 
The stationary distribution of this conditioned random variable is the 
so-called quasi-stationary distribution  (QSD). 
The QSD is important for the information it gives about the long-term 
behavior of a surviving population. 
Our main goal is to derive approximations of the first three cumulants of 
the QSD.

Our approach toward this goal resembles the well-known moment closure method 
in the sense that both are based on a system of ODEs for the first few 
(three) cumulants. 
This system of ODEs is not closed: 
The number of unknown functions (cumulants) is larger than the number of 
equations. 
This fact has caused all practitioners of moment closure to introduce an 
approximating assumption that leads to closure of the system of ODEs.   
The method is used widely, as shown by the recent review paper by Kuehn (2016).
However, it is also known to possess several weaknesses. 
The first of these weaknesses is that no condition for validity of the 
resulting approximation is available, a second one is that the magnitude of 
the approximation error can not be evaluated, a third one is that it often 
leads to spurious solutions that require large efforts to eliminate, and a 
fourth one is that the dependence of the approximating expressions on the 
model's parameters is not known. 
An early basic paper in the area was written by Whittle (1957). 
The desirability of closing the system of ODEs for the cumulants cannot be 
denied. 
It is indeed a necessity if the goal is to derive exact expressions for the 
cumulants.
However, it is important to realize that closure is not necessary if one is 
satisfied with approximate results, as we are in this paper.

We do not use moment closure methods at all. 
Instead, we apply the alternative mathematical method introduced in 
N{\aa}sell (2017). 
Thus, we base our work on a system of ODEs for cumulants, just as is done by 
practitioners of moment closure. 
However, we avoid the step that leads to closure of the system of ODEs, since 
the ad-hoc nature of this step is the cause of several weaknesses of the 
moment closure method. 
Our alternative is to search for asymptotic approximations. 
This requires the reparametrization that accompanies the second one of the two 
model formulations of Section 2. 
To proceed with our approach, we also need information about the orders of 
magnitude (in terms of $N$) of the cumulants that appear in the system of ODEs. 
We motivate our assumptions in this regard with numerical evaluations. 
The method that we use avoids the weaknesses associated with moment 
closure methods.

Methods for numerical evaluations of the QSD and of the associated cumulants 
are discussed in Section 3. 
Our approach is different from what has been used by previous workers on 
this model. 
Results of numerical evaluations are used to motivate forms of basic 
assumptions that are made for determining asymptotic approximations of the 
first few cumulants of the QSD, and also for checking the orders of 
magnitude of the errors of the approximations that we derive.

Sections 4 and 5 are used to study the two random variables $X(t)$ and 
$X^Q(t)$, respectively.  
We give a system of first-order ODEs for the first 3 cumulants of the 
unconditioned random variable $X(t)$ in Section 4.
The same system of equations is shown in Section 5 to provide an approximation 
for the system of ODEs of the first 3 cumulants of the conditioned random 
variable $X^Q(t)$ in a particular parameter region. 
This means that approximations of the cumulants of the QSD in this region are found 
as critical points of the system of ODEs. 
Asymptotic approximations of the coordinates of the critical points of this 
system of ODEs are derived in Section 6. 
They serve to give asymptotic approximations of the first 3 cumulants of the QSD. 
These results hold for small positive integer values of the parameter $s$ that 
describes the power of the population size $n$ that gives the decreasing 
function of the population size that reflects the density dependence of 
the net birth rate of the model.   
An extension of these results to arbitrary positive values of $s$, both 
integer and non-integer, is proposed in Section 7. 
Comparisons are made in Section 8 with published results that also deal with 
the cumulants of the QSD of the same model. 
The paper ends with some concluding comments in Section 9.

\section{Model Formulation}

The stochastic power law logistic model is a model for the size of a 
population with density-dependent growth. 
It is formulated as a birth-death process $\{X(t),t\ge 0\}$. 
The hypotheses of the model are summarized in the descriptions of the 
population birth-rate $\lambda_n$ and the population death-rate $\mu_n$ as 
functions of the state $n$ of the process.
The formulation given by Matis, Kiffe, and Parthasarathy (1998) takes the 
following form:
The population birth-rate equals  
\begin{equation} \label{2.1} 
   \lambda_n = \begin{cases}                                                    
                            (a_1 - b_1 n^s) n,   &  n\le (a_1/b_1)^{1/s},  \\
                                         0,             & \text{otherwise},
                       \end{cases}
\end{equation}
and the population death-rate is 
\begin{equation} \label{2.2}
     \mu_n = (a_2 + b_2 n^s)n.
\end{equation}
The state space appears to be unbounded and equal to $\{0,1,2,\dots\}$. 
The parameters $a_1, b_1, a_2, b_2, s$ are all assumed to be positive. 
This model formulation takes its form from the similar model introduced by 
Bartlett (1957) for the special case of the Verhulst logistic model where 
$s=1$.

We shall use a different formulation, where 
\begin{align} \label{2.3}
    & \lambda_n = \mu R_0 \left(1 - \left(\frac{n}{N}\right)^s \right)n, 
    \quad n=0,1,\dots,N, \\ \label{2.4}
    & \mu_n       = \mu \left(1 + \alpha \left(\frac{n}{N}\right)^s \right)n, 
\quad n=0,1,\dots,N.
\end{align}
The state space of the process in this formulation is finite and equal to 
$\{0,1,2,\dots,N\}$. 
The parameter space contains the five parameters $N, R_0, \alpha, \mu$, 
and $s$. 
Among these, $N$ is a large positive integer that represents the maximum 
population size, $R_0$ is a positive dimensionless threshold parameter, 
$\alpha$ is a nonnegative dimensionless parameter, $\mu$  is a positive 
death rate with dimension inverse time, and $s$ is a positive dimensionless 
number.   
This model formulation is an extension from the formulation of the logistic 
Verhulst model, corresponding to $s=1$, given by N{\aa}sell (2017). 
Under the assumption that the initial distribution is supported on the state 
space, the process remains there for all time, since $\mu_0=0$ and 
$\lambda_N=0$. 
The origin is seen to be an absorbing state, since both $\lambda_0$ and 
$\mu_0$ are equal to zero. 
Absorption at the origin corresponds to extinction of the population. 
To study a surviving population, we introduce the conditioned random variable 
$X^{(Q)}(t)$ by conditioning $X(t)$ on the event $X(t)>0$. 
Thus, 
\begin{equation} \label{2.5}
    P\{X^{(Q)}(t)=k\} =  P\{X(t)=k|X(t)>0\}, \quad k=1,2,\dots,N.
\end{equation}
The state space of the conditioned random variable $X^{(Q)}(t)$ differs from that of 
$X(t)$ in one respect: The origin can be reached by $X(t)$, but not by the 
conditioned random variable $X^{(Q)}(t)$. 
Thus, the state space of the latter one of these random variables is equal to $\{1,2,\dots,N\}$. 
Stationary distributions of the two random variables are vastly different.  
The stationary distribution of $X(t)$ is degenerate with probability one at 
the origin, while the stationary distribution of the conditioned random 
variable $X^{(Q)}(t)$ is the important quasi-stationary distribution (QSD), 
supported on the state space $\{1,2,\dots,N\}$.

We shall use the second of the two model formulations in this paper. 
An important reason for this is that the second formulation contains a parameter $N$ 
that can take large values. 
The presence of such a parameter is essential for the formulation of asymptotic approximations. 
The second formulation also gives knowledge about the parameter dimensions. 
Dimensions of parameters and state variables need to be known when different terms 
are compared.  
It is obvious that only terms of the same dimension are comparable.

The two formulations of the stochastic power law logistic model are essentially equivalent 
if we require in the first formulation that $(a_1/b_1)^{1/s}=N$ is an integer, that the state 
space is finite and equal to $\{0,1,2,\dots,N\}$, and that $b_2 \ge 0$. 
We assume in what follows in this paper that these three requirements are met. 
The second model formulation can then be seen as a reparametrization of the first one. 
We find then that the 4 parameters $a_1, a_2, b_1, b_2$ of the first formulation can be 
expressed in terms of the 5 parameters $\mu, R_0, \alpha, N, s$ of the second 
formulation as follows:
\begin{align} \label{2.6}
   & a_1 =  \mu R_0, \\ \label{2.7}
   & a_2 =  \mu, \\ \label{2.8}
   & b_1 =  \mu \frac{R_0}{N^s}, \\ \label{2.9} 
   & b_2 =  \mu \frac{\alpha}{N^s}.
\end{align}
Obviously, $s$ takes the same value in both formulations. 
We note also that some results that we shall refer to below use the following 
notations for the sums and differences of the parameters $a_i$ and $b_i$: 
\begin{align} \label{2.10} 
   & a = a_1 - a_2 = \mu(R_0-1), \\ \label{2.11}
   & b = b_1 + b_2 = \mu\frac{R_0+\alpha}{N^s}, \\ \label{2.12}
   & c = a_1 + a_2 = \mu(R_0+1), \\ \label{2.13}
   & d = b_1 - b_2 = \mu \frac{R_0-\alpha}{N^s}. 
\end{align}

As already mentioned, the logistic Verhulst model can be seen as a special 
case of the power law logistic model, corresponding to $s=1$. 
Asymptotic approximations of the first 3 cumulants of the QSD for the 
Verhulst model are derived in N{\aa}sell (2017). 
The goal of the present paper is to derive corresponding approximations of 
the first 3 cumulants of the QSD of the power law logistic model for 
positive values of $s$. 
We proceed by first deriving asymptotic approximations of the first 3 
cumulants of the QSD for the $s$-values 2, 3, and 4.  
These results are then extended in Section 7 to both integer and non-integer 
positive values of $s$.

\section{Numerical Evaluations}

Numerical evaluations have been used in earlier work on this model, and we 
shall also use them here. 
However, there are differences between our approach and those of earlier 
workers both in the methods used for deriving numerical results, and in the 
type of evaluations that are carried out.

Bartlett {\it et al.} (1960) argue incorrectly that a stationary distribution 
cannot exist for the model that they are concerned with, since $\lambda_0=0$, 
while we claim that the process in this case has a degenerate stationary 
distribution with probability one at the origin. 
Bartlett {\it et al.} claim furthermore that "the probability distribution 
for the stationary (or quasi-stationary) distribution must satisfy the 
recurrence relation"
\begin{equation}  \label{3.1}
   \mu_n p_n = \lambda_{n-1} p_{n-1}.
\end{equation}
We find it disturbing that they refer to both stationary and quasi-stationary distributions here, 
without clarifying which of these two distributions that they mean. 
The statement can certainly not hold for both of them. 
The fact is that it is incorrect for both of them. 
Instead, the relation \eqref{3.1} holds for the stationary distribution 
$p^{(0)}$ of a related auxiliary process $X^{(0)}(t)$ that is useful for 
studying the quasi-stationary distribution. 
The transition rates of this auxiliary process are equal to those of the 
process $X(t)$, with the one exception that the rate $\mu_1$ for transition 
from the state 1 to the absorbing state 0 is replaced by zero.  
The stationary distribution $p^{(0)}$ is non-degenerate and can be evaluated 
explicitly. 
It has been used to study the QSD for a SIS model, which is the special 
case of the model that we are concerned with here that corresponds to the 
parameter values $s=1$ and $\alpha=0$, in N{\aa}sell (2011).  
It follows from this study that $p^{(0)}$ is a good approximation of the 
QSD in its body in the parameter region where $R_0>1$. 
But it follows also that $p^{(0)}$ is not an acceptable approximation of the QSD if $R_0<1$, nor in the 
parameter region near $R_0=1$. 
Similar relations between $p^{(0)}$ and the QSD are expected to hold for 
$\alpha>0$ and for  integer values of $s$ larger than 1. 
Both Bartlett {\it et al.} (1960) and Matis {\it et al.} (1998) use the 
recurrence relation \eqref{3.1} as a basis for numerical evaluations 
of the QSD. 
Our approach is different, as shown below. 
However, we expect the cumulants of the stationary distribution 
$p^{(0)}$ to be close to the cumulants of the QSD in the parameter region where 
$R_0>1$, since the tail probablities have a minor influence on the 
cumulant values. 
We note that Renshaw (2011) uses  $p^{(0)}$ as a "working definition" 
of quasi-stationarity. 
We do not follow this approach. 
Our results require a sharp line to be drawn between exact and 
approximate relations.

To describe our method for numerical evaluations, we need a second 
auxiliary process $X^{(1)}(t)$.  
The birth rates of this process are equal to those of the original 
process $X(t)$, while the death rates are slightly smaller than those of 
$X(t)$, and equal to $\mu_{n-1}$. 
This means that the second auxiliary process can be interpreted as the 
original process $X(t)$ with one surviving and immortal individual.  
Like the first of the two auxiliary processes, the second one lacks absorbing 
state, and its stationary distribution $p^{(1)}$ is non-degenerate and can be 
evaluated explicitly.

Our numerical evaluations of the QSD make use of the restart map $\Psi$ 
analyzed by Ferrari, Kesten, Martinez, and Picco (1995), and also discussed 
by N{\aa}sell (2011). 
This restart map is defined as follows: 
Let a probability vector $\nu = (\nu_1, \nu_2, \dots, \nu_N)$ be given. 
Define a process related to the process $X(t)$ that we are studying by the 
requirement that whenever the original process reaches the state zero, it is 
immediately restarted at some state $j$ with probability $\nu_j$. 
The restarted process has the state space $\{1, 2, \dots, N\}$, and a unique 
stationary distribution $p$. 
The map $\Psi$ is then defined by $\Psi(\nu)=p$. 
The quasi-stationary distribution $q$ is a fixed point of the mapping $\Psi$, 
since $\Psi(q)=q$. 
Furthermore, Ferrari {\it et al.} (1995) show that an iteration scheme can 
be established by repeated applications of the map $\Psi$ to an arbitrary 
initial probability vector, and that it converges to the quasi-stationary 
distribution $q$. 
Suitable initial probability vectors are given by the stationary 
distributions $p^{(0)}$ and $p^{(1)}$  of the two auxiliary processes 
$X^{(0)}(t)$ and $X^{(1)}(t)$.   
It turns out that that the map $\Psi$ can be described explicitly with 
the aid of these two stationary distributions. 
A derivation is given by N{\aa}sell (2011). 
By denoting the components of the vector $p=\Psi(\nu)$ by $p_n$, we find that 
\begin{equation} \label{3.2} 
   p_n = \pi_n S_n p_1, \quad n=1,2,\dots,N, 
\end{equation}
where 
\begin{equation} \label{3.3}
   \pi_1=1, \quad \pi_n = \frac{\lambda_1\lambda_2\cdots\lambda_{n-1}}    
{\mu_2\mu_3\cdots\mu_n}, \quad n=2,3,\dots,N,
\end{equation}
and 
\begin{equation} \label{3.4}
   S_n = \sum_{k=1}^{n} \frac{1-\sum_{j=1}^{k-1}\nu_j}{\rho_k}, \quad n=1,2,\dots,N,
\end{equation}
with
\begin{equation} \label{3.5}
   \rho_1=1, \quad \rho_n 
= \frac{\lambda_1\lambda_2\cdots\lambda_{n-1}}{\mu_1\mu_2\cdots\mu_{n-1}}, 
\quad n=2,3,\dots,N,
\end{equation}
and
\begin{equation} \label{3.6}
   p_1 = \frac{1}{\sum_{n=1}^N \pi_n S_n}.
\end{equation}
The stationary distributions $p^{(0)}$ and $p^{(1)}$ of the two auxiliary processes are 
determined from tte sequencies $\pi_n$ and $\rho_n$ as follows:  
\begin{equation} \label{3.7} 
    p_n^{(0)} = \pi_n p_1^{(0)}, \quad n=1,2,\dots,N, \quad \text{where} \quad
        p_1^{(0)} = \frac{1}{\sum_{n=1}^N \pi_n},
\end{equation}
and 
\begin{equation} \label{3.8} 
    p_n^{(1)} = \rho_n p_1^{(1)}, \quad n=1,2,\dots,N, \quad \text{where} \quad
        p_1^{(1)} = \frac{1}{\sum_{n=1}^N \rho_n},
\end{equation}

Our numerical method for determining the quasi-stationary distribution 
consists in applying the restart map $\Psi$ to a suitable initial 
distribution and continue iterations until successive iterates are 
sufficiently close. 
In case $R_0>1$, which is the parameter region of main interest for the 
study of cumulants in this paper, the stationary distribution $p^{(0)}$ 
is recommended over the distribution $p^{(1)}$ as initial distribution.

Numerical evaluations are used by Matis {\it et al.} (1998) for showing 
that the errors that their approximations lead to are small. 
With our different parametrization we have access to a parameter $N$ that 
can take large values. 
We can then use numerical evaluations to extract additional information 
both about the model and about the approximations that we derive. 
One goal is to derive information about the orders of magnitude (in terms 
of $N$) of various cumulants, and another one concerns the probability of 
taking the value 1 in the quasi-stationary distribution. 
Information of these two types is important for the formulation of 
hypotheses as a basis for finding asymptotic approximations of the first 
few cumulants for large values of $N$, as will be shown below. 
Another use of numerical evaluations is to derive information about the orders 
of magnitude (in terms of $N$) of the errors that are committed in using 
our approximations. 
Such results will be used below to support the forms of asymptotic 
approximations that we derive.

\section{ODEs for Cumulants of the Unconditioned Random Variable $X(t)$}

The starting point for our derivations of asymptotic approximations of the first 3 
cumulants of the QSD are ODEs for the first 3 cumulants of the conditioned random
variable $X^Q(t)$.  
It turns out that these ODEs are closely related to the ODEs for the corresponding 
cumulants of the unconditioned random variable $X(t)$ in one important region of 
parameter space, namely where $R_0>1$. 
We start therefore by giving these latter ODEs.

The derivation of these ODEs makes use of the moment generating function $M(\theta,t)$ and the  
cumulant generating function $K(\theta,t)$ of the random variable $X(t)$. 
They are defined by 
\begin{equation} \label{4.1}
   M(\theta,t) = E \exp(\theta X(t))  = \sum_{n=0}^N p_n(t) \exp(n\theta)
\end{equation}
and 
\begin{equation} \label{4.2} 
   K(\theta,t) = \log M(\theta,t),
\end{equation}
where the state probabilities $p_n(t)=P\{X(t)=n\}$ satisfy the Kolmogorov forward equations 
\begin{equation} \label{4.3}
    p_n'(t) = \lambda_{n-1} p_{n-1}(t) - (\lambda_n + \mu_n) p_n(t) + \mu_{n+1} p_{n+1}(t), 
                \quad n = 0,1,\dots,N.
\end{equation} 
(Here, we agree to put $\lambda_{-1} = \mu_{N+1} = p_{-1}(t) = p_{N+1}(t) = 0$ so that \eqref{4.3} 
makes sense formally for all n-values indicated.)

The definitions of the cumulant generating function $K$ and of the transition rates $\lambda_ n$ and 
$\mu_n$ can be used to derive a partial differential equation (PDE) for $K$. 
It can be written as follows: 
\begin{multline} \label{4.4}
     \frac{\partial K(\theta,t)}{\partial t} = \mu (\exp(\theta)-1) 
        \left[(R_0-\exp(-\theta)) \frac{\partial K(\theta,t)}{\partial \theta} \right. \\
           \left. - \frac{R_0 + \alpha \exp(-\theta)}{N^s} \exp(-K(\theta,t)) 
             \frac{\partial^{s+1}\exp(K(\theta,t))}{\partial\theta^{s+1}} \right].
\end{multline} 
The same result follows by applying the reparametrization in \eqref{2.6}--\eqref{2.9} to the expression for the PDE given by Renshaw (2011), and numbered (3.4.11).

The cumulants $\kappa_i(t)$ are obtained from a power series expansion of the cumulant generating function $K(\theta,t) $ as follows:  
\begin{equation} \label{4.5} 
   K(\theta,t) = \sum_{i=1}^{\infty} \frac{\kappa_i(t)}{i!} \theta^i.
\end{equation} 

We use this result for determining the ODEs of the first few cumulants of the random variable $X(t)$. 
They are found by determining the first few terms of the expansion of the PDE of the cumulant generating function $K(\theta,t)$ in terms of $\theta$, and then identifying terms with equal powers of $\theta$.

To express our results, we use the capital letters $A$, $B$, and $C$ to denote the time derivatives of the first three cumulants of the random variable $X(t)$, with superscripts indicating the $s$-values, as follows: 
\begin{align} \label{4.6}
    A^{(s)}(t) & = \kappa_1^{(s)^{'}}(t), \\ \label{4.7}
    B^{(s)}(t) & = \kappa_2^{(s)^{'}}(t), \\ \label{4.8}
    C^{(s)}(t) & = \kappa_3^{(s)^{'}}(t).
\end{align}

Expressions for the functions $A^{(s)}, B^{(s)}, C^{(s)}$ for the $s$-values 2, 
3 and 4 can then be written as follows:

\begin{align} \label{4.9}
   A^{(2)} & =  \mu(R_0-1) \kappa_1 - \mu\frac{R_0 + \alpha}{N^2} 
     (\kappa_1^3 + 3 \kappa_1 \kappa_2 + \kappa_3), \\ \label{4.10}
   B^{(2)} & =  2 \mu(R_0-1) \kappa_2 +  \mu (R_0+1)\kappa_1  
       - \mu\frac{R_0 - \alpha}{N^2} (\kappa_1^3 + 3 \kappa_1 \kappa_2 
         + \kappa_3)  \\* \notag
      & \phantom{as} - 2 \mu\frac{R_0 + \alpha}{N^2} 
       (3 \kappa_1^2 \kappa_2 + 3 \kappa_1 \kappa_3 + 3 \kappa_2^2 + \kappa_4),  
     \\ \label{4.11}
  C^{(2)} & = \mu (R_0-1) (\kappa_1 + 3 \kappa_3) + 3 \mu (R_0+1) \kappa_2 
     \\* \notag
     & \phantom{as} - 3 \mu \frac{R_0 - \alpha}{N^2} (3 \kappa_1^2 \kappa_2 
        + 3 \kappa_1 \kappa_3 + 3 \kappa_2^2 + \kappa_4)   \\* \notag 
      &  \phantom{as} - \mu \frac{R_0 + \alpha}{N^2} (\kappa_1^3 
     + 9 \kappa_1^2 \kappa_3 + 18 \kappa_1 \kappa_2^2 + 3 \kappa_1 \kappa_2   
     + 9 \kappa_1 \kappa_4   + 27 \kappa_2 \kappa_3    \\* \notag
         & \phantom{hejahejaheja}  + \kappa_3 + 3 \kappa_5), \\ \label{4.12}
   A^{(3)} & =  \mu (R_0-1) \kappa_1  
         - \mu\frac{R_0 + \alpha}{N^3}  (\kappa_1^4 + 6 \kappa_1^2 \kappa_2 
           + 4 \kappa_1 \kappa_3 + 3 \kappa_2^2 + \kappa_4 ), \\ \label{4.13}
   B^{(3)} & =   2 \mu (R_0-1) \kappa_2 + \mu (R_0+1)\kappa_1  \\* \notag
      & \phantom{as} - \mu \frac{R_0 - \alpha}{N^3} (\kappa_1^4 
     + 6 \kappa_1^2 \kappa_2 + 4 \kappa_1 \kappa_3 + 3 \kappa_2^2 + \kappa_4)   
    \\* \notag
      & \phantom{as} - 2 \mu \frac{R_0 + \alpha}{N^3} 
       (4 \kappa_1^3 \kappa_2 + 6 \kappa_1^2 \kappa_3  
          + 12 \kappa_1 \kappa_2^2   + 4 \kappa_1 \kappa_4  
         + 10 \kappa_2 \kappa_3 
        + \kappa_5),  \\ \label{4.14}
  C^{(3)} & =   \mu(R_0-1) (\kappa_1 + 3 \kappa_3) + 3 \mu(R_0+1) \kappa_2 
   \\* \notag
     &  \phantom{as} - 3 \mu \frac{R_0 - \alpha}{N^3} (4 \kappa_1^3 \kappa_2 
        + 6 \kappa_1^2 \kappa_3 + 12 \kappa_1 \kappa_2^2  + 4 \kappa_1 \kappa_4  
           + 10 \kappa_2 \kappa_3  + \kappa_5)   \\* \notag 
      &  \phantom{as} -  \mu \frac{R_0 + \alpha}{N^3} (\kappa_1^4 
      + 12 \kappa_1^3 \kappa_3  + 36 \kappa_1^2 \kappa_2^2 
      + 6 \kappa_1^2 \kappa_2 + 18 \kappa_1^2 \kappa_4 
       + 108 \kappa_1 \kappa_2 \kappa_3  \\*\notag
            & \phantom{hejahejaheja} 
             + 4 \kappa_1 \kappa_3 + 12 \kappa_1 \kappa_5   
          + 36 \kappa_2^3 + 3 \kappa_2^2 + 42 \kappa_2 \kappa_4  + 30 \kappa_3^2  
\\*\notag
             & \phantom{hejahejaheja} + \kappa_4 + 3 \kappa_6 ), \\ \label{4.15}
   A^{(4)} & =  \mu (R_0-1) \kappa_1  \\* \notag
      & \phantom{as} - \mu\frac{R_0+\alpha}{N^4}(\kappa_1^5 
     + 10 \kappa_1^3 \kappa_2  + 10 \kappa_1^2 \kappa_3 
     + 15 \kappa_1 \kappa_2^2 
   + 5 \kappa_1 \kappa_4 + 10 \kappa_2 \kappa_3 + \kappa_5),   \\ \label{4.16}
   B^{(4)} & =  2 \mu(R_0-1)\kappa_2 +  \mu(R_0+1) \kappa_1  \\*\notag
   & \phantom{he} - \mu\frac{R_0-\alpha}{N^4}(\kappa_1^5 + 10\kappa_1^3 \kappa_2 
      + 10 \kappa_1^2 \kappa_3 + 15 \kappa_1 \kappa_2^2
     + 5 \kappa_1 \kappa_4 + 10 \kappa_2 \kappa_3 + \kappa_5) \\*\notag
      & \phantom{he} - 2 \mu\frac{R_0+\alpha}{N^4}(5\kappa_1^4 \kappa_2 
     + 10 \kappa_1^3 \kappa_3 + 30 \kappa_1^2 \kappa_2^2 + 10 \kappa_1^2 \kappa_4 
     + 50 \kappa_1 \kappa_2 \kappa_3 + 5 \kappa_1 \kappa_5 \\*\notag
       & \phantom{hejhejtjejer} + 15 \kappa_2^3  + 15 \kappa_2 \kappa_4 
      + 10 \kappa_3^2 + \kappa_6), \\ \label{4.17}
   C^{(4)} & =  \mu (R_0-1)(\kappa_1 + 3 \kappa_3) + 3 \mu (R_0+1) \kappa_2 
    \\*\notag
    & \phantom{as} - 3  \mu\frac{R_0-\alpha}{N^4} (5 \kappa_1^4 \kappa_2 
    + 10 \kappa_1^3 \kappa_3 + 30 \kappa_1^2 \kappa_2^2 + 10 \kappa_1^2 \kappa_4 
       + 50 \kappa_1 \kappa_2 \kappa_3 + 5 \kappa_1 \kappa_5 \\*\notag
       & \phantom{hejakillarman}+ 15 \kappa_2^3  + 15 \kappa_2 \kappa_4 
     + 10 \kappa_3^2 + \kappa_6) \\*\notag
      & \phantom{he} -  \mu\frac{R_0+ \alpha}{N^4} 
      (\kappa_1^5 + 15 \kappa_1^4 \kappa_3 + 60 \kappa_1^3 \kappa_2^2 
     + 10 \kappa_1^3 \kappa_2 + 30 \kappa_1^3 \kappa_4 
           + 270 \kappa_1^2 \kappa_2 \kappa_3 \\*\notag
      & \phantom{hejatjejern}+ 10 \kappa_1^2 \kappa_3 + 30 \kappa_1^2 \kappa_5 
          + 180 \kappa_1 \kappa_2^3 + 15 \kappa_1 \kappa_2^2 
            + 210 \kappa_1 \kappa_2 \kappa_4 \\*\notag
      & \phantom{hejatjejern}+ 150 \kappa_1 \kappa_3^2 + 5 \kappa_1 \kappa_4 
        + 15 \kappa_1 \kappa_6  + 285 \kappa_2^2 \kappa_3 + 10 \kappa_2 \kappa_3 
             + 60 \kappa_2 \kappa_5 \\*\notag
      & \phantom{hejatjejern}+ 105 \kappa_3 \kappa_4 + \kappa_5 + 3 \kappa_7).
\end{align}

Derivations of these results using Maple are given in N{\aa}sell (2018). 

We note that the 3 derivatives $A^{(2)}, A^{(3)}, A^{(4)}$ in \eqref{4.9}, \eqref{4.12}, \eqref{4.15} 
of the first cumulant $\kappa_1$ for the 3 $s$-values 2, 3, and 4 are found from the 
expressions (27), (28), (29) in Matis, Kiffe, Parthasarathy (1998) by reparametrization 
using \eqref{2.10}--\eqref{2.13} after introducing the correction that the first term 
$(a-b\kappa_1^4)$ in the right-hand side of (29) is written $(a-b\kappa_1^4)\kappa_1$. 
Similarly, the 2 derivatives $B^{(2)}$ and $B^{(3)}$ in \eqref{4.10} and \eqref{4.13} 
of the second cumulant $\kappa_2$ for the 2 $s$-values 2 and 3 are found from the 
expressions (33) and (34) in Matis {\it et al.} by the same reparametrization after 
noting that the minus sign of the term $-2a\kappa_2$ in the right-hand side of their 
formula (34) is incorrect and should be changed to a plus sign.

\section{ODEs for Cumulants of the Conditioned Random Variable $X^{(Q)}(t)$}

The goal of our study is to determine approximations of the first three cumulants of the QSD 
for the stochastic variable $X(t)$. 
This leads us to a study of the stationary values of the system of ODEs for the first three 
cumulants of the conditioned random variable $X^{(Q)}(t)$. 
These ODEs turn out to be closely related to the ODEs for the cumulants of the unconditioned 
random variable $X(t)$ in one important parameter region, namely where $R_0>1$. 
Identificaton of this parameter region is of high importance, since it gives the 
validity condition for the approximations that we derive.  

In similarity to the case in the previous section we use the cumulant generating function 
$K^{(Q)}(\theta,t)$  of the random variable $X^{(Q)}(t)$ that is of concern here to derive 
the ODEs of the first few cumulants of this random variable. 
To derive the PDE of this cumulant generating function we proceed as in the previous section. 
The main difference from the case in the previous section is that the expression for the time derivative 
of the state probability $p_n^{(Q)}(t)$ is different from the counterpart of \eqref{4.3}. 
By using the relation in \eqref{2.5} we find that
\begin{equation} \label{5.1}
       p_n^{(Q)}(t) = \frac{p_n(t)}{1-p_0(t)}, \quad n=1,2,\dots,N.
\end{equation}
 Differentiation and use of the Kolmogorov forward equations in \eqref{4.3} gives 
\begin{multline} \label{5.2}
       {p_n^{(Q)}}'(t) =    \lambda_{n-1} p_{n-1}^{(Q)}(t) - (\lambda_n + \mu_n) p_n^{(Q)}(t) 
        + \mu_{n+1} p_{n+1}^{(Q)}(t) \\ 
         + \mu_1 p_1^{(Q)}(t) p_n^{(Q)}(t),  \quad n=1,2,\dots,N. 
\end{multline}

The PDE for the cumulant generating function $K^{(Q)}(\theta,t)$ of the 
conditioned random variable $X^Q(t)$ turns out to be quite similar to the 
PDE for $K(\theta,t)$.  
It can be written as follows: 
\begin{multline} \label{5.3}
     \frac{\partial K^{(Q)}(\theta,t)}{\partial t} = \mu (\exp(\theta)-1) 
       \left[(R_0-\exp(-\theta)) \frac{\partial K^{(Q)}(\theta,t)}{\partial \theta}
        \right. \\
           \left. - \frac{R_0 + \alpha \exp(-\theta)}{N^s} \exp(-K^{(Q)}(\theta,t)) 
                  \frac{\partial^{s+1}\exp(K^{(Q)}(\theta,t))}{\partial\theta^{s+1}} 
           \right] \\
         + \mu\left(1+\frac{\alpha}{N}\right) p_1^{(Q)}(t) 
          \left(1 - \exp(-K^{(Q)}(\theta,t)\right).
\end{multline}

As in the previous section we determine the ODEs of the first few cumulants by expanding this PDE 
in terms of $\theta$ and identifying terms of equal powers of $\theta$. 
We find then that the last term of the right-hand side of every such ODE 
contains $p^{(Q)}_1(t)$ as a factor. 
The stationary values of the probabilities $p_n^{(Q)}(t)$ are equal to the 
quasi-stationary probabilities $q_n$. 
We claim that the probability $q_1$ is exponentially small in $N$ for any $\alpha\ge0$ and 
any positive integer value of $s$ when $R_0>1$. 
Arguments that support this are given below. 
The last term can therefore be ignored when we are searching for asymptotic approximations of 
the critical points.  
An important consequence of this is that asymptotic approximations of the 
first few cumulants of the QSD for $R_0>1$ are found as stationary solutions 
of the ODEs in Section 4 for the corresponding cumulants of the unconditioned 
random variable $X(t)$, with obvious rules for exclusion of spurious solutions. 
Derivations of these results are given below in Section 6.

To show that $q_1$ is exponentially small in $N$ when $R_0>1$ for any $\alpha\ge0$ 
and any positive integer $s$ we note first from N{\aa}sell (2011) that this result is true for 
the SIS model, where $\alpha=0$ and $s=1$. 
Next we consider the Verhulst model, with $\alpha>0$ and $s=1$. 
We then make use of the result that $q_1<p_1^{(0)}$, where $p_n^{(0)}$ is the stationary 
distribution of the auxiliary process $X^{(0)}(t)$, and where it is shown by N{\aa}sell (2001b) that 
\begin{equation} \label{5.4} 
     p_1^{(0)} \sim \frac{R_0(R_0-1)\sqrt{(1+\alpha)N}}{R_0+\alpha}\varphi(\beta_1), 
     \quad R_0>1, \quad N \to \infty,
\end{equation}
where
\begin{equation} \label{5.5} 
    \beta_1 = \sqrt{2N \gamma_1}
\end{equation}
and 
\begin{equation} \label{5.6} 
     \gamma_1 = \log R_0 - \frac{\alpha+1}{\alpha} \log \frac{(\alpha+1)R_0}{R_0+\alpha}, 
      \quad \alpha>0. 
\end{equation}
Since the normal density function $\varphi$ satisfies $\varphi(y) = \exp(-y^2/2)/\sqrt{2\pi}$, 
we conclude that $p_1^{(0)}$, and therefore also $q_1$, are exponentially small in $N$ 
under the condition $R_0>1$, $s=1$ and $\alpha>0$.

Finally, we consider the cases where $s$ is an arbitrary positive integer, $\alpha\ge 0$, and $R_0>1$. 
Here, we have not found an analytical proof that $q_1$ is exponentially small. 
Instead, we use numerical evaluations to support our claim. 
The results in case $R_0=2$, $\alpha=1$, the four $s$-values 1,2,3,4, and the three $N$-values 100, 200, and 400, are given in Table 1. 
They show that $q_1$ decreases strongly with $s$ for the indicated values of $R_0$, $\alpha$, and $N$.  
Since $q_1$ is exponentially small in $N$ for $s=1$, we conclude that $q_1$ is exponentially small in 
$N$ for the indicated parameter values.   
The results in Table 1 also support our claim in the sense that a doubling of the $N$-value for fixed 
$s$ leads to an approximate squaring of the value of $q_1$.  
This is what is expected when $q_1$ is exponentially small in $N$. 

The condition that we have found under which $q_1$ is exponentially small in $N$, namely $R_0>1$, turns out to be the condition for the validity of the approximations of the first three cumulants of the QSD that we derive later. 
We note that a corresponding condition of validity has been missing from earlier results for the same model, based on cumulant closure.

\begin{table}[h]
  \begin{center}
    \begin{tabular}{  | c | r | r | r | }
       \hline
            $s$  & $N=100$ \phantom{A} & $N=200$ \phantom{A}  & $N=400$ \phantom{A}
\\ \hline
    1     & $1.30*10^{-5}$       & $1.49*10^{-10}$      & $1.27*10^{-20}$  \\
    2     & $6.59*10^{-12}$     & $18.2*10^{-24}$      & $95.7*10^{-48}$  \\
    3     & $7.74*10^{-16}$     & $18.5*10^{-32}$      & $73.7*10^{-64}$ \\ 
    4     & $2.18*10^{-18}$     & $1.22*10^{-36}$      & $0.264*10^{-72}$  \\
 \hline
\end{tabular}
   \vskip 4mm
    \caption{Numerical evaluations of the probability $q_1$  of the QSD of the 
                 stochastic power law logistic model. 
                 Results are shown for $R_0=2$, $\alpha=1$, the $s$-values 
                 1, 2, 3, and 4, and the $N$-values 100, 200, and 400. 
                 The results indicate that $q_1$ is exponentially small. 
                 Derivations of these results using Maple are given in N{\aa}sell (2018).}
  \end{center}
\end{table}

\section{Derivations of Asymptotic Approximations of Cumulants of the QSD} 

This section is used to derive asymptotic approximations of the first 3 
cumulants of the quasi-stationary distribution of the stochastic power 
law logistic model for the $s$-values 2, 3, and 4 and large $N$-values 
in the parameter region where $R_0>1$. 
The method that we use is similar to the one introduced in the study of 
the cumulants of the QSD of the Verhulst logistic model in N{\aa}sell (2017). 
This latter study actually represents the special case with $s=1$ of the 
model that we study here. 
Our results are based on assumptions about the forms of asymptotic 
approximations of the first few cumulants for large values of $N$. 
For the Verhulst model with $s=1$ and $R_0>1$ we give arguments in 
N{\aa}sell (2017) for assuming that the first 4 cumulants of the QSD 
are $\text{O(N)}$. 
The results make it easy to conjecture that all finite order cumulants of 
the QSD for $s=1$ are $\text{O}(N)$ in the parameter region $R_0>1$. 
It is tempting to guess that this holds also for integer values of $s$ 
larger than 1.

\begin{table}[h]
  \begin{center}
    \begin{tabular}{ | c | c | r | r | r | }
       \hline
    $s$ &  Cumulant  & N=100 & N=200 & N=400\\ \hline
     1 & $\kappa_1$  & 81.6    & 163     & 327 \\
     1 & $\kappa_2$  & 16.7    & 33.3    & 66.3 \\
     1 & $\kappa_3$  &-13.8    & -27.3   & -54.3 \\ 
     1 & $\kappa_4$  & 8.61    & 16.9    & 33.6  \\
     1 & $\kappa_5$  & -0.532 & -0.620 & -0.833 \\
     1 & $\kappa_6$  &-10.0   & -21.0   & -42.8 \\  
     1 & $\kappa_7$  & 17.8   & 38.3    & 79.0 \\ \hline
     4 & $\kappa_1$  & 95.0   & 190     & 380 \\
     4 & $\kappa_2$  & 4.93   & 9.74    & 19.3  \\
     4 & $\kappa_3$  & -4.80  & -9.45   & -18.8 \\
     4 & $\kappa_4$  & 4.63   & 9.09    & 18.0 \\ 
     4 & $\kappa_5$  & -4.61  & -8.98  & -17.7 \\ 
     4 & $\kappa_6$  & 5.52   & 10.5   &  20.5 \\ 
     4 & $\kappa_7$  &  -10.2 & -18.9 & -36.6 \\ \hline
    \end{tabular}
    \vskip 4mm
    \caption{Numerical evaluations of the first 7 cumulants of the QSD of 
     the stochastic power law logistic model. 
     Results are shown for $R_0=10$, $\alpha=1$, the $s$-values 1 and 4, 
     and the $N$-values 100, 200, and 400. 
     The results indicate that the first 7 cumulants are all $\text{O}(N)$.  
     }
  \end{center}
\end{table}

In what follows in this section we shall assume that the first $s+3$ 
cumulants of the QSD are $\text{O(N)}$ when $s=1, 2, 3, 4$, $\alpha\ge0$,  
and $R_0>1$. 
As a basis for this we refer to the numerical results in Table 2. 
It shows numerically determined values of the first 7 cumulants of the QSD 
for the parameter values $R_0=10$ and $\alpha=1$, with the 2 $s$-values 1 
and 4, and with the three $N$-values 100, 200, and 400. 
It is easy to verify from the table that all cumulants experience an 
approximate doubling both when $N$ is doubled from 100 to 200 and also 
when $N$ is doubled from 200 to 400. 
We take this as a strong indication that all of the cumulants considered 
are $\text{O}(N)$.  
The results in Table 2 are derived using Maple in N{\aa}sell (2018), 
where also similar results are given for the $s$-values 2 and 3.

We study first the case $s=2$.  
The assumptions that the first 5 cumulants $\kappa_1-\kappa_5$ are 
$\text{O}(N)$ in the parameter region $R_0>1$ for $s =2$ are basic for 
our further assumptions that the first 5  cumulants have the following 
asymptotic behaviors for $R_0>1$ and $s=2$:

\begin{align} \label{6.1}
  & \kappa_1 = x_1 N + x_2 + \frac{x_3}{N} + \text{O}\left(\frac{1}{N^2}\right),
     \quad R_0>1, \quad s=2, \\ \label{6.2} 
  & \kappa_2 = y_1 N + y_2 + \text{O}\left(\frac{1}{N}\right), \quad R_0>1, 
  \quad s=2, \\    \label{6.3}
   & \kappa_3 = z_1 N + \text{O}(1), \quad R_0>1, \quad s=2, \\ \label{6.4}
   & \kappa_4 = u_1 N + \text{O}(1), \quad R_0>1, \quad s=2, \\ \label{6.5}
   & \kappa_5 = u_2 N + \text{O}(1), \quad R_0>1, \quad s=2. 
\end{align}
The reason for including different numbers of terms in these asymptotic 
approximations will become apparent shortly. 
Thus, we have introduced 8 unknowns, namely $x_1, x_2, x_3, y_1, y_2, z_1, u_1, 
u_2$. 
We proceed to determine the first 6 of them. 
By inserting the above asymptotic approximations of the first 5 cumulants into 
the expressions 
\eqref{4.4}--\eqref{4.6} for the functions $A^{(2)}, B^{(2)}, C^{(2)}$, 
we find that asymptotic approximations for them can be written as follows:
\begin{align}  \label{6.6}
   & A^{(2)} = A^{(2)}_1 N + A^{(2)}_2 + \frac{A^{(2)}_3}{N} 
+ \text{O}\left(\frac{1}{N^2}\right), 
    \quad R_0>1, \quad s=2, \\ \label{6.7}
   & B^{(2)} = B^{(2)}_1 N + B^{(2)}_2 + \text{O}\left(\frac{1}{N}\right),   
\quad R_0>1,  \quad s=2, \\ \label{6.8}
   & C^{(2)} = C^{(2)}_1 N + \text{O}(1), \quad R_0>1, \quad s=2,  
\end{align}
where 
\begin{align} \label{6.9}
   &  A^{(2)}_1 = \mu (R_0-1) x_1 - \mu (R_0+\alpha) x_1^3, \\  \label{6.10}
 &  A^{(2)}_2 =\mu (R_0-1) x_2 - 3\mu(R_0+\alpha) (x_1^2 x_2 + x_1 y_1), 
\\  \label{6.11}
   &  A^{(2)}_3 = \mu(R_0-1) x_3 -\mu(R_0+\alpha) 
(3 x_1^2 x_3 + 3 x_1 x_2^2 + 3 x_1 y_2 
                           + 3 x_2 y_1 + z_1), \\   \label{6.12}
   & B^{(2)}_1 = 2 \mu(R_0-1) y_1 +\mu (R_0+1) x_1  -\mu (R_0-\alpha) x_1^3
               - 6\mu (R_0 + \alpha) x_1^2 y_1,     \\  \label{6.13}
 & B^{(2)}_2 = 2 \mu(R_0-1) y_2  + \mu(R_0+1) x_2 - 3\mu(R_0-\alpha) 
(x_1^2 x_2 + x_1 y_1) 
\\ \notag
 & \phantom{abcdef}   - 6 \mu(R_0+\alpha) (x_1^2 y_2 + 2 x_1 x_2 y_1 + x_1 z_1 
+ y_1^2), \\ \label{6.14}
  & C^{(2)}_1  = \mu(R_0-1) (x_1 + 3 z_1) + 3 \mu(R_0+1) y_1 - 9 \mu(R_0-\alpha) 
x_1^2 y_1
                    \\  \notag
   & \phantom{abcdef}  - \mu(R_0+\alpha) (x_1^3 + 9 x_1^2 z_1 + 18 x_1 y_1^2).   
\end{align}
The reason for including different numbers of terms in the assumed asymptotic 
approximations of the cumulants $\kappa_1 - \kappa_5$ can be understood with 
reference to the resulting asymptotic approximations of the quantities 
$A^{(2)}, B^{(2)}$, and $C^{(2)}$. 
We note for example from the expression \eqref{4.4} for $A^{(2)}$ that an 
additional term in the assumed asymptotic approximation of $\kappa_2$ would 
contribute a term of the order of $1/N^2$ to the asymptotic approximation of 
$A^{(2)}$. 
Inclusion of such an additional term would therefore be absorbed in the error 
term for the asymptotic approximation of $A^{(2)}$. 
We note from the expressions in \eqref{6.9}--\eqref{6.14} that only the 
first 6 of the 8 coefficients introduced above appear in these expressions.

Our assumptions for $s=2$ that the first 5 cumulants have the asymptotic approximations 
given in \eqref{6.1}--\eqref{6.5} lead to the asymptotic approximations of $A^{(2)}$, 
$B^{(2)}$, $C^{(2)}$ in \eqref{6.6}--\eqref{6.8}. 
These expressions are equal to the derivatives with respect to time of the first 3 cumulants. 
Setting them equal to zero gives conditions for the stationary points of the corresponding 
system of ODEs. 
Since we are working with approximations instead of exact results, we transform these 
conditions into the 3 requirements that $A^{(2)}=\text{O}(1/N^2)$, 
$B^{(2)}=\text{O}(1/N)$, $C^{(2)}=\text{O}(1)$. 
These requirements are satisfied by the basic conditions that the 6 expressions 
$A^{(2)}_1, A^{(2)}_2, A^{(2)}_3, B^{(2)}_1, B^{(2)}_2, C^{(2)}_1$ in 
\eqref{6.9}--\eqref{6.14} are equal to zero. 
We note that each of these expressions is a polynomial of degree 3 in 
the 6 unknown coefficients $x_1, x_2, x_3, y_1, y_2, z_1$. 
Our basic problem is to solve the 6 equations that are formed by setting the corresponding 
expressions equal to zero for the 6 unknown coefficients.   
Each equation is a polynomial of degree $s+1=3$ on the 6 unknown coefficients. 
Solution appears possible, since the number of equations equals the number of unknowns. 
Further inspection of the 6 equations reveals that considerable simplification can be achieved in the 
solution by determining the 6 unknowns in a definite order, as described below. 
The first advantage is that each of the equations to be solved contains only one unknown coefficient, 
and the second one is that all equations except the first one are linear in the coefficient to be solved for. 
The very first equation is of degree $s+1=3$ in the case we are considering here, with $s=2$.
All of the roots of the first equation except one are spurious solutions that turn out to be easy to identify and to delete. 
Whenever a solution is found for one of the coefficients, its value is immediately inserted into the remaining unsolved equations. 
The order of solution that leads to these pleasant results is as follows: 
First, the equation $A^{(2)}_1=0$ is solved for $x_1$. 
This equation has in this case, with $s=2$, $s+1=3$ solutions. 
Among them, we exclude $x_1=0$ and $x_1<0$ as the only spurious solutions that appear in this method. 
The solution found for $x_1$ is then inserted into the remaining 5 unsolved equations.  
After finding $x_1$, we proceed to solve the equation $B^{(2)}_1=0$ for $y_1$, 
the equation $A^{(2)}_2=0$ for $x_2$, the equation $C^{(2)}_1=0$ for $z_1$, 
the equation $B^{(2)}_2=0$ for $y_2$, and the equation $A^{(2)}_3=0$ for $x_3$. 
It is elementary to use these rules to solve for the 6 unknown coefficients. 
Alternative derivations using Maple are given in N{\aa}sell (2018).  
The solutions are as follows:

\begin{align} \label{6.15}
   & x_1 = \left(\frac{R_0-1}{R_0+\alpha}\right) ^{1/2}, 
      \quad R_0>1, \quad \alpha\ge0, \quad s=2, \\  \label{6.16}
   & x_2 = -\frac34 \frac{(\alpha+1)R_0}{(R_0+\alpha)(R_0-1)}, 
      \quad R_0>1, \quad \alpha\ge0, \quad s=2, \\  \label{6.17}
   & x_3 = - \frac{(\alpha+1)R_0}{(R_0+\alpha)^2 (R_0-1)^2}
                  \left( \frac{R_0+\alpha}{R_0-1} \right)^{1/2}     \\   \notag
   &  \phantom{hejsanhej} \cdot \left[\frac78 
  (R_0^2+\alpha)+\frac{13}{32}(\alpha+1)R_0\right], 
       \quad R_0>1, \quad \alpha\ge0,   \quad s=2,\\  \label{6.18}
   & y_1 = \frac12 \frac{(\alpha+1)R_0}{(R_0+\alpha) (R_0-1)} 
            \left(\frac{R_0-1}{R_0+\alpha}\right)^{1/2}, 
       \quad R_0>1, \quad \alpha\ge0,  \quad s=2,\\  \label{6.19}
 & y_2 = \frac34 \frac{(\alpha+1)R_0} {(R_0+\alpha)^2 (R_0-1)^2} (R_0^2+\alpha),  
        \quad R_0>1, \quad \alpha\ge0,  \quad s=2, \\  \label{6.20}
   & z_1 = - \frac{(\alpha+1)R_0} {(R_0+\alpha)^2(R_0-1)^2} 
                \left(\frac{R_0-1}{R_0+\alpha} \right)^{1/2}    \\  \notag    
 &  \phantom{hejsanhej}    \cdot  \left[\frac12(R_0^2+\alpha)-\frac14(\alpha+1)
  R_0\right],  
         \quad R_0>1, \quad \alpha\ge0, \quad  s=2. 
\end{align}

The method used here for derivation of asymptotic approximations of the first 
3 cumulants of the QSD for $s=2$ can in principle be followed for higher 
integer values of $s$. 
We proceed to consider the case when $s=3$. 
We begin by asssuming that the first 6 cumulants have the following asymptotic 
behaviors for $R_0>1$ and $s=3$:

\begin{align} \label{6.21}
  & \kappa_1 = x_1 N + x_2 + \frac{x_3}{N} + \text{O}\left(\frac{1}{N^2}\right), 
                   \quad R_0>1, \quad s=3, \\ \label{6.22}
   & \kappa_2 = y_1 N + y_2 + \text{O}\left(\frac{1}{N}\right), 
     \quad R_0>1, \quad s=3, \\     \label{6.23}
   & \kappa_3 = z_1 N + \text{O}(1), \quad R_0>1, \quad s=3, \\ \label{6.24}
   & \kappa_4 = u_1 N + \text{O}(1), \quad R_0>1, \quad s=3, \\ \label{6.25}
   & \kappa_5 = u_2 N + \text{O}(1), \quad R_0>1, \quad s=3,  \\  \label{6.26}
   & \kappa_6 = u_3 N  + \text{O}(1), \quad R_0>1, \quad s=3.
\end{align}

Insertions of these asymptotic approximations of the first 6 cumulants into 
the expressions \eqref{4.7}-\eqref{4.9} for the functions $A^{(3)}, B^{(3)}, 
C^{(3)}$ lead to the following asymptotic approximations for them:
\begin{align}  \label{6.27}
   & A^{(3)} = A^{(3)}_1 N + A^{(3)}_2 + \frac{A^{(3)}_3}{N} 
+ \text{O}\left(\frac{1}{N^2}\right),  \quad R_0>1, \quad s=3, \\ \label{6.28}
   & B^{(3)} = B^{(3)}_1 N + B^{(3)}_2 + \text{O}\left(\frac{1}{N}\right), 
    \quad R_0>1, \quad s=3, \\  \label{6.29}
   & C^{(3)} = C^{(3)}_1 N + \text{O}(1), \quad R_0>1, \quad s=3,  
\end{align}
where
\begin{align} \label{6.30}
   &  A^{(3)}_1 = \mu(R_0-1) x_1 - \mu(R_0+\alpha) x_1^4, \\  \label{6.31}
   &  A^{(3)}_2 = \mu(R_0-1) x_2 - 2\mu(R_0+\alpha) (2x_1^3 x_2 +3x_1^2 y_1), \\  
\label{6.32}
 & A^{(3)}_3 = \mu(R_0-1) x_3 -\mu(R_0+\alpha) 
(4 x_1^3 x_3 + 6 x_1^2 x_2^2 + 6 x_1^2 y_2  + 12 x_1 x_2 y_1     \\  \notag 
&  \phantom{hejsan} + 4 x_1 z_1+ 3 y_1^2), \\   \label{6.33}
   & B^{(3)}_1 = 2 \mu(R_0-1) y_1 + \mu(R_0+1) x_1  - \mu(R_0-\alpha) x_1^4  
                       - 8 \mu(R_0 + \alpha) x_1^3 y_1,\\  \label{6.34}
& B^{(3)}_2 = 2 \mu(R_0-1) y_2 + \mu(R_0+1) x_2 - \mu(R_0-\alpha) 
   (4 x_1^3 x_2 + 6 x_1^2 y_1) \\ \notag
&  \phantom{abcdef} - \mu(R_0+\alpha) 
(8 x_1^3 y_2 + 24 x_1^2 x_2 y_1 + 12 x_1^2 z_1 + 24 x_1 y_1^2),   \\ \label{6.35}  
   & C^{(3)}_1 = \mu(R_0-1) (x_1 + 3 z_1) + 3 \mu(R_0+1) y_1 
  - 12 \mu(R_0-\alpha) x_1^3 y_1 \\ \notag
& \phantom{abcdef} - \mu(R_0+\alpha) (x_1^4 + 12 x_1^3 z_1 + 36 x_1^2 y_1^2).   
\end{align}

These 6 expressions are all polynomials of degree 4 in the 6 unknowns  
$x_1$, $x_2$, $x_3$, $y_1$, $y_2$, $z_1$.
We solve the 6 equations formed by putting each of these expressions equal 
to zero for the 6 unknowns $x_1, x_2, x_3, y_1, y_2, z_1$. 
First, the equation $A^{(3)}_1=0$ is solved for $x_1$. 
Among the 4 roots we exclude the one that is equal to zero, and the two that 
are complex conjugate as the only spurious solutions that appear in this 
method. 
The remaining 5 equations are then solved for the remaining 5 unknowns, 
using an order of solution similar to what is described for the case $s=2$ 
above.  
Derivations using Maple are given in N{\aa}sell (2018).  
The results are as follows: 
\begin{align} \label{6.36}
   & x_1 = \left(\frac{R_0-1}{R_0+\alpha}\right)^{1/3},  
     \quad R_0>1, \quad \alpha\ge 0, \quad s=3, \\  \label{6.37}
   & x_2 = -\frac23 \frac{(\alpha+1)R_0}{(R_0+\alpha)(R_0-1)}, 
      \quad R_0>1, \quad \alpha\ge 0, \quad s=3,\\  \label{6.38}
   & x_3 = - \frac{(\alpha+1)R_0}{(R_0+\alpha)^2(R_0-1)^2} 
              \left(\frac{R_0+\alpha}{R_0-1} \right)^{1/3}   \\  \notag
    &   \phantom{hejsanhej}  \cdot \left[\frac89(R_0^2+\alpha)
+\frac{7}{27}(\alpha+1)R_0\right], 
       \quad R_0>1, \quad \alpha\ge0,    \quad s=3\\ \label{6.39}
   & y_1 = \frac13 \frac{(\alpha+1)R_0}{(R_0+\alpha) (R_0-1)}
                 \left(\frac{R_0-1}{R_0+\alpha} \right)^{1/3}, 
      \quad R_0>1, \quad \alpha\ge0,  \quad s=3,  \\  \label{6.40}
   & y_2 = \frac23 \frac{(\alpha+1)R_0}{(R_0+\alpha)^2 (R_0-1)^2}(R_0^2+\alpha), 
       \quad R_0>1,  \quad \alpha\ge0,     \quad s=3, \\  \label{6.41}
   & z_1 = - \frac{(\alpha+1)R_0}{(R_0+\alpha)^2(R_0-1)^2}
                  \left(\frac{R_0-1}{R_0+\alpha}\right)^{1/3}    \\  \notag
   & \cdot \phantom{hejsanhej} \cdot \left[\frac13(R_0^2+\alpha)
-\frac19(\alpha+1)R_0\right],  
      \quad R_0>1, \quad \alpha\ge0, \quad s=3. 
\end{align}  

The last problem addressed in this section is the derivation of asymptotic 
approximations of the first 3 cumulants of the QSD for $s=4$. 
We begin by assuming that the first 7 cumulants have the following 
asymptotic behaviors for $R_0>1$ and $s=4$:

\begin{align} \label{6.42}
  & \kappa_1 = x_1 N + x_2 + \frac{x_3}{N} + \text{O}\left(\frac{1}{N^2}\right), 
                   \quad R_0>1, \quad s=4, \\ \label{6.43}
   & \kappa_2 = y_1 N + y_2 + \text{O}\left(\frac{1}{N}\right), 
     \quad R_0>1, \quad s=4, \\     \label{6.44}
   & \kappa_3 = z_1 N + \text{O}(1),  \quad R_0>1, \quad s=4, \\ \label{6.45}
   & \kappa_4 = u_1 N + \text{O}(1),  \quad R_0>1, \quad s=4, \\ \label{6.46}
   & \kappa_5 = u_2 N + \text{O}(1),  \quad R_0>1, \quad s=4,  \\  \label{6.47}
   & \kappa_6 = u_3 N  + \text{O}(1), \quad R_0>1, \quad s=4, \\ \label{6.48}
   & \kappa_7 = u_4 N + \text{O}(1),  \quad R_0>1, \quad s=4. 
\end{align}

Insertions of these asymptotic approximations of the first 7 cumulants into 
the expressions \eqref{4.10}--\eqref{4.12}  for the functions $A^{(4)}, 
B^{(4)}, C^{(4)}$ lead to the following asymptotic approximations for these 
3 functions:

\begin{align}  \label{6.49}
   & A^{(4)} = A^{(4)}_1 N + A^{(4)}_2 + \frac{A^{(4)}_3}{N} 
    + \text{O}\left(\frac{1}{N^2}\right),     \quad R_0>1, \quad s=4, \\ 
\label{6.50}
   & B^{(4)} = B^{(4)}_1 N + B^{(4)}_2 + \text{O}\left(\frac{1}{N}\right), 
    \quad R_0>1, \quad s=4, \\      \label{6.51}
   & C^{(4)} = C^{(4)}_1 N + \text{O}(1), \quad R_0>1, \quad s=4,  
\end{align}
where 
\begin{align} \label{6.52}
   A^{(4)}_1 & = \mu(R_0-1) x_1 - \mu(R_0+\alpha) x_1^5, \\  \label{6.53} 
   A^{(4)}_2 & = \mu(R_0-1) x_2 - 5 \mu(R_0+\alpha) (x_1^4 x_2 + 2 x_1^3 y_1), \\ 
\label{6.54}
   A^{(4)}_3 & = \mu(R_0-1) x_3 - 5 \mu(R_0+\alpha) 
   (x_1^4 x_3 + 2 x_1^3 x_2^2 + 2 x_1^3 y_2  + 6 x_1^2 x_2 y_1 \\*\notag
     & \phantom{he}  + 2 x_1^2 z_1 + 3 x_1 y_1^2), \\ \label{6.55}
   B^{(4)}_1 & = 2\mu(R_0-1)y_1  + \mu(R_0+1) x_1 - \mu(R_0-\alpha) x_1^5
       -10 \mu(R_0+\alpha) x_1^4 y_1,   \\ \label{6.56}
   B^{(4)}_2 & = 2 \mu(R_0-1) y_2 + \mu(R_0+1) x_2 - 5 \mu(R_0-\alpha) 
   (x_1^4 x_2 + 2 x_1^3 y_1) \\ \notag
 & \phantom{he} - 10 \mu(R_0+\alpha)  (x_1^4 y_2 + 4 x_1^3 x_2 y_1 + 2 x_1^3 z_1 
       + 6 x_1^2 y_1^2),   \\ \label{6.57}
   C^{(4)}_1 & = \mu(R_0-1) (x_1 + 3 z_1) + 3 \mu(R_0+1) y_1 - 15 \mu(R_0-\alpha) 
  x_1^4 y_1 \\ \notag
  & \phantom{he}  - \mu(R_0+\alpha) (x_1^5 + 15 x_1^4 z_1 + 60 x_1^3 y_1^2).
 \end{align}

These 6 expressions are all polynomials of degree 5 in the 6 unknown quantities 
$x_1, x_2, x_3, y_1, y_2, z_1$. 
As above, we solve the 6 equations formed by putting each of these expressions 
equal to zero for the 6 unknowns $x_1, x_2, x_3, y_1, y_2, z_1$. 
First, the equation $A^{(4)}_1 = 0$ is solved for $x_1$. 
Among the 5 roots we exclude 4 of them as spurious solutions, namely one that 
is equal to zero, one that is negative, and 2 that are imaginary. 
The remaining 5 equations are then solved for the remaining 5 unknowns, 
using the same order as above for the cases $s=2$ and $s=3$. 
Derivations using Maple are given by N{\aa}sell (2018). 
The results are as follows:
\begin{align} \label{6.58} 
   x_1 = & \left(\frac{R_0-1}{R_0+\alpha}\right)^{1/4},  
     \quad R_0>1, \quad \alpha\ge0, \quad s=4, \\ \label{6.59}
   x_2 = & - \frac{5}{8} \frac{(\alpha+1)R_0}{(R_0+\alpha)(R_0-1)}, 
     \quad R_0>1, \quad \alpha\ge0, \quad s=4,  \\ \label{6.60}
   x_3 = & - \frac{(\alpha+1)R_0}{(R_0+\alpha)^2 (R_0-1)^2} 
              \left(\frac{R_0+\alpha}{R_0-1} \right)^{1/4}   \\  \notag
& \phantom{hejsan}  
\cdot \left[ \frac{15}{16}(R_0^2+\alpha) + \frac{25}{128}(\alpha+1)R_0 \right], 
    \quad R_0>1, \quad \alpha\ge0,\quad s=4, \\ \label{6.61}
   y_1 = & \frac14 \frac{(\alpha+1)R_0}{(R_0+\alpha) (R_0-1)} 
          \left(\frac{R_0-1}{R_0+\alpha}\right)^{1/4}, 
    \quad R_0>1, \quad \alpha\ge0,\quad s=4, \\ \label{6.62}
  y_2 = & \frac58 \frac{(\alpha+1)R_0}{(R_0+\alpha)^2(R_0-1)^2} (R_0^2+\alpha), 
      \quad R_0>1, \quad \alpha\ge0,    \quad s=4, \\ \label{6.63}
   z_1 = & - \frac{(\alpha+1)R_0}{(R_0+\alpha)^2(R_0-1)^2}
             \left(\frac{R_0-1}{R_0+\alpha} \right)^{1/4}  \\   \notag
  &  \phantom{hejsan}   \cdot \left[ \frac14 (R_0^2+\alpha) 
    - \frac{1}{16}(\alpha+1)R_0 \right], 
       \quad R_0>1, \quad \alpha\ge0,   \quad s=4.  
\end{align}

\section{Extensions to Positive Values of $s$}
 
The results derived in the previous section give approximations of the first 3 cumulants of the QSD 
for the integer $s$-values 2, 3, and 4. 
The present section is used to extend these approximations to positive values of $s$. 
The extension is made in four steps. 
In the very first step we combine the results of the previous section with the results in 
N{\aa}sell (2017), which are valid for $s=1$, to give results that hold for the $s$-values 1, 2, 3, and 4. 
In the second step we extend these results to the integer $s$-values from 1 to 10. 
The third step extends these results to all positive integer values of $s$.  
The final step makes an extension from this to all positive values of $s$.

We emphasize that the method that we use for deriving approximations  of the first few cumulants 
for integer values of $s$ cannot be used for deriving results for non-integer values of $s$. 
The reason is that the derivations are based on ODEs for the first few cumulants of integer order. 
Our last step of extension does not require us to consider ODEs for cumulants of non-integer order. 
Instead we use a continuity argument for extending the derived approximations of the first few 
cumulants to non-integer positive values of $s$.

By using obvious definitions of the 6 coefficients 
$x_1, x_2, x_3, y_1, y_2, z_1$, we find from N{\aa}sell (2017) for the case 
$s=1$ that they depend on the model parameters $R_0$ and $\alpha$ as follows: 

\begin{align} \label{7.1}
   & x_1 = \frac{R_0-1}{R_0+\alpha}, 
    \quad R_0>1, \quad \alpha\ge0, \quad s=1, \\  \label{7.2}
   & x_2 = - \frac{(\alpha+1)R_0}{(R_0+\alpha)(R_0-1)},  
    \quad R_0>1, \quad \alpha\ge0,  \quad s=1,\\  \label{7.3}
   & x_3 = - \frac{(\alpha+1)R_0} {(R_0-1)^3} (R_0+1),
    \quad R_0>1, \quad \alpha\ge0,   \quad s=1,\\  \label{7.4}
   & y_1 = \frac{(\alpha+1)R_0}{(R_0+\alpha)^2}, 
    \quad R_0>1, \quad \alpha\ge0,   \quad s=1, \\  \label{7.5}
   & y_2 =  \frac{(\alpha+1)R_0} {(R_0+\alpha)^2 (R_0-1)^2} (R_0^2+\alpha),  
    \quad R_0>1, \quad \alpha\ge0,   \quad s=1,  \\  \label{7.6}
   & z_1 = -  \frac{(\alpha+1)R_0} {(R_0+\alpha)^3} (R_0-\alpha) ,
    \quad R_0>1, \quad \alpha\ge0,   \quad s=1. 
\end{align}

The dependence on $s$ for all positive values of $s$ is claimed to be as follows for the 6 coefficients 
$x_1, x_2, x_3, y_1, y_2, z_1$: 
\begin{align} \label{7.7}
   & x_1 = \left( \frac{R_0-1}{R_0+\alpha} \right)^{h_1(s)},
    \quad R_0>1, \quad \alpha\ge0,  \quad s > 0, \\  \label{7.8}
   & x_2 = - h_2(s) \frac{(\alpha+1) R_0}{(R_0+\alpha)(R_0-1)},
    \quad R_0>1, \quad \alpha\ge0,  \quad s > 0,  \\  \label{7.9}
   & x_3 = - \frac{(\alpha+1)R_0}{(R_0+\alpha)^2(R_0-1)^2}              
                \left(\frac{R_0+\alpha}{R_0-1} \right)^{h_1(s)}   \\  \notag
 &  \phantom{hejsanhej} \cdot 
     \left[h_3(s) (R_0^2+\alpha)+h_4(s)(\alpha+1)R_0\right], \\ \notag
  &  \phantom{hejhejsan} \quad R_0>1, \quad \alpha\ge0,  \quad s > 0,\\ \label{7.10}
   & y_1 = h_1(s) \frac{(\alpha+1)R_0}{(R_0+\alpha)(R_0-1)} 
                \left(\frac{R_0-1}{R_0+\alpha} \right)^{h_1(s)},
    \quad R_0>1, \quad \alpha\ge0,  \quad s > 0, \\  \label{7.11}
   & y_2 = h_2(s) \frac{(\alpha+1)R_0}{(R_0+\alpha)^2(R_0-1)^2} 
     (R_0^2+\alpha), \quad R_0>1, \quad \alpha\ge0,  \quad s > 0,  \\  
\label{7.12}
   & z_1 = -\frac{(\alpha+1)R_0}{(R_0+\alpha)^2(R_0-1)^2}  
                  \left(\frac{R_0-1}{R_0+\alpha} \right)^{h_1(s)}  \\   \notag
 &  \phantom{hejsanhej} 
    \cdot \left[ h_1(s) (R_0^2+\alpha) - h_5(s)(\alpha+1)R_0 \right], 
     \quad R_0>1, \quad \alpha\ge0,  \quad s > 0, 
\end{align}

where the functions $h_1 - h_5$ are defined as follows:
\begin{align}  \label{7.13}
   h_1(s) = & \frac{1}{s},  \\  \label{7.14}
   h_2(s) = & \frac{s+1}{2s},  \\   \label{7.15}
   h_3(s) = & \frac{s^2 + 6s + 5}{12 s},  \\  \label{7.16}
   h_4(s) = & \frac{s^2 + 12s +11}{24 s^2},  \\  \label{7.17}
   h_5(s) = & \frac{1}{s^2}.
\end{align}

To show that the functions $h_1-h_5$ are determined by these expressions for $s > 0$, we determine first the values that they take for the integer $s$-values 1 - 10, and given in Table 3. 
The values taken by these 5 functions for the integer $s$-values 1-4 are readily found from the expressions  found for the 6 coefficients  $x_1, x_2, x_3, y_1, y_2, z_1$ given in the previous and present sections.  
The further results for the $s$-values from 5 to 10 are based on asymptotic approximations of the first 
3 cumulants of the QSD of the power law logistic model that have been derived using Maple, and are
reported in N{\aa}sell (2018). 

It is straightforward to use the entries in Table 3 for the functions $h_1$, $h_2$, and $h_5$ to confirm that the expressions \eqref{7.13}, \eqref{7.14}, \eqref{7.17} for these functions are valid for all integer
values of $s$ from 1 to 10.

\begin{table}[t]
  \begin{center}
    \begin{tabular}{ | c | c | c | c | c | c |}
       \hline
    $s$ &  $h_1(s)$  & $h_2(s)$ & $h_3(s)$ & $h_4(s)$ & $h_5(s)$ \\ \hline
     1   &       1        &    1         &    1         &     1        &  1           \\
     2   &       1/2     &    3/4      &    7/8      &  13/32     &  1/4       \\
     3   &       1/3     &    2/3      &    8/9      &    7/27     &  1/9       \\  
     4   &       1/4     &    5/8      &   15/16   & 25/128     &  1/16     \\
     5   &       1/5     &    3/5      &   1         &  4/25        &  1/25     \\
     6   &       1/6     &    7/12    &   77/72   & 119/864    &  1/36    \\
     7   &       1/7     &    4/7      &   8/7      & 6/49          &  1/49    \\
     8   &       1/8     &    9/16    &   39/32   & 57/512      &   1/64   \\
     9   &       1/9     &    5/9      &   35/27   & 25/243      &   1/81   \\
    10  &      1/10    &    11/20   &   11/8     & 77/800      &   1/100 \\ \hline
      \end{tabular}
    \caption{Values of the functions $h_1 - h_5$  when their argument 
     $s$ takes integer values from 1 to 10. 
     Derivations for $s$ from 5 to 10 are given in N{\aa}sell (2018).  }
  \end{center}
\end{table}

We turn now to deal with the functions $h_3$ and $h_4$. 
To find expressions for them, we assume that each of them can be written as a quotient of 2 
polynomials in $s$ of degree 2. 
 To be specific, we assume that 
\begin{align} \label{7.18}
   h_3(s) = \frac{a_0 + a_1 s + a_2 s^2}{b_0 + b_1 s + b_2 s^2}, \\ \label{7.19}
   h_4(s) = \frac{c_0 + c_1 s + c_2 s^2}{d_0 + d_1 s + d_2 s^2}. 
\end{align}
To determine the values of  the 6 coefficients $a_0, a_1, a_2, b_0, b_1, b_2$, we establish 6 equations 
that are found by equating the values of $h_3(s)$ for the integer $s$-values from 1 to 6 according to \eqref{7.18} to the corresponding values in Table 3. 
By solving these equations using Maple, as seen in N{\aa}sell (2018), we find that $h_3(s)$ equals 
the expression in \eqref{7.15}. 
An entirely similar treatment of the function $h_4$ leads to the expression \eqref{7.16} for 
$h_4(s)$. 
So far, we can conclude that the expressions \eqref{7.15} and \eqref{7.16} for $h_3(s)$ and $h_4(s)$ 
are valid for the integer $s$-values in the interval from 1 to 6.  
Extensions of the validities of the 2 expressions for $h_3(s)$ and $h_4(s)$ to all integer values of $s$ 
from 1 to 10 are easily seen to hold by the simple expedient of verifying equalities between the 
values of the expressions in \eqref{7.15} and \eqref{7.16} and the corresponding values in Table 3. 
We take this as a strong indication that the domain of validity of the expressions \eqref{7.15} and \eqref{7.16} can be further extended to all positive integer values of $s$.

As a last step of extension, we conjecture that the domain of validity for the 5 functions $h_1 - h_5$ can 
be extended from all positive integers $s$ to all positive values of $s$. 
We illustrate the result for non-integer $s$-values by giving the numerical values of the error terms for the 3 cumulants $\kappa_1, \kappa_2, \kappa_3$ in case $R_0=10$ and $\alpha=1$, for the 2 $s$-values 
0.5 and 3.5, and the 3 $N$-values 100, 200, and 400 in Table 4. 
The table indicates that the error term for $\kappa_1$ is divided by approximately 4 for each doubling 
of $N$, the error term for $\kappa_2$ is divided by approximately 2 for each doubling of $N$, and 
the error term for $\kappa_3$ is approximately constant when $N$ is doubled. 
This is consistent with the claims that the error terms of $\kappa_k$ are $\text{O}(1/N^{3-k})$ 
for the $k$-values 1, 2, and 3. 
Evaluations of the error terms in Table 4 are derived using Maple in N{\aa}sell (2018).

\begin{table}[t]
  \begin{center}
    \begin{tabular}{ | c | c | r | r | r | }
       \hline
     $s$ &  Cumulant  & N=100 \phantom{ab} & N=200 \phantom{ab} & N=400 
\phantom{ab} \\ \hline
     0.5 & $\kappa_1$   & $-227*10^{-6}$   & $-55*10^{-6}$  & $-14*10^{-6}$ \\
     0.5 & $\kappa_2$   & $111*10^{-4}$    & $54*10^{-4}$   & $27*10^{-4}$ \\
     0.5 & $\kappa_3$   & $-35*10^{-2}$    & $-34*10^{-2}$  & $-33*10^{-2}$ \\  
\hline
     3.5 & $\kappa_1$   & $-334*10^{-7}$  & $-83*10^{-7}$   & $-21*10^{-7}$ \\
     3.5 & $\kappa_2$   & $247*10^{-5}$   & $122*10^{-5}$  & $61*10^{-5}$ \\
     3.5 & $\kappa_3$   & $-144*10^{-3}$  & $-142*10^{-3}$ & $-141*10^{-3}$ \\ 
\hline
      \end{tabular}
    \vskip 4mm
    \caption{Numerical evaluations of the error terms of the approximations 
of the first 3 cumulants of the QSD of the stochastic power law logistic model. 
     Results are shown for $R_0=10$, $\alpha=1$, the $s$-values 0.5 and 3.5,  
and the $N$-values 100, 200, and 400.     
The table indicates that the error term of $\kappa_1$ is $\text{O}(1/N^2)$, 
the error term of $\kappa_2$ is  $\text{O}(1/N)$, and the error term of 
$\kappa_3$ is $\text{O}(1)$. }
  \end{center}
\end{table}

The results in \eqref{7.7}--\eqref{7.12} show the way in which our asymptotic approximations 
of the first 3 cumulants of the QSD of the stochastic power law logistic model depend 
on the 4 parameters $N$, $R_0$, $\alpha$, and $s$. 
The one-term asymptotic approximations of the first 3 cumulants 
$\kappa_1$, $\kappa_2$, and $\kappa_3$ are useful for the information that they give about the 
behaviors of these cumulants. 
They can be written as follows: 
\begin{align} \label{7.20}
   \kappa_1 = x_1 N + O(1), \quad R_0>1, \quad \alpha\ge 0, \quad s > 0, \quad N\to\infty, \\ \label{7.21}
   \kappa_2 = y_1 N + O(1), \quad R_0>1, \quad \alpha\ge 0, \quad s > 0, \quad N\to\infty, \\ \label{7.22}
   \kappa_3 = z_1 N + O(1), \quad R_0>1, \quad \alpha\ge 0, \quad s > 0, \quad N\to\infty, 
\end{align}
where expressions for the 3 coefficients $x_1, y_1, z_1$ are given in \eqref{7.7}, \eqref{7.10}, and \eqref{7.12}, respectively. 
We use these expressions to derive some properties of the first 3 cumulants, 
valid for sufficiently large values of $N$. 

We note first that the expectation $\kappa_1$ is an increasing function of the power $s$. 
This follows from the following expression for the derivative of $x_1$ with respect to $s$:
\begin{equation} \label{7.23}
    \frac{dx_1}{ds} = \frac{1}{s^2} \log \left( \frac{R_0+\alpha}{R_0-1} \right) 
      \left( \frac{R_0-1}{R_0+\alpha} \right)^{1/s}.
\end{equation}

Our second observation concerns the variance $\kappa_2$. 
For sufficiently large values of $N$, we find that it has a maximum as a function of $s$ at $s=s_2$, 
where 
\begin{equation} \label{7.24}
     s_2 = \log ((R_0+\alpha)/(R_0-1)). 
\end{equation}
This follows from the following expression for the derivative of $y_1$ with respect to $s$:
\begin{equation} \label{7.25}
   \frac{dy_1}{ds} = - \frac{1}{s^3} \frac{(\alpha+1)R_0}{(R_0+\alpha)(R_0-1)}
     \left( \frac{R_0-1}{R_0+\alpha} \right)^{1/s} \left(s - \log \frac{R_0+\alpha}{R_0-1} \right).
\end{equation}
We conclude in particular that $\kappa_2$ is increasing in $s$ for $s < s_2$ and decreasing in $s$
for $s>s_2$, provided $N$ is large enugh.

\begin{figure}[h]
\epsfig{file=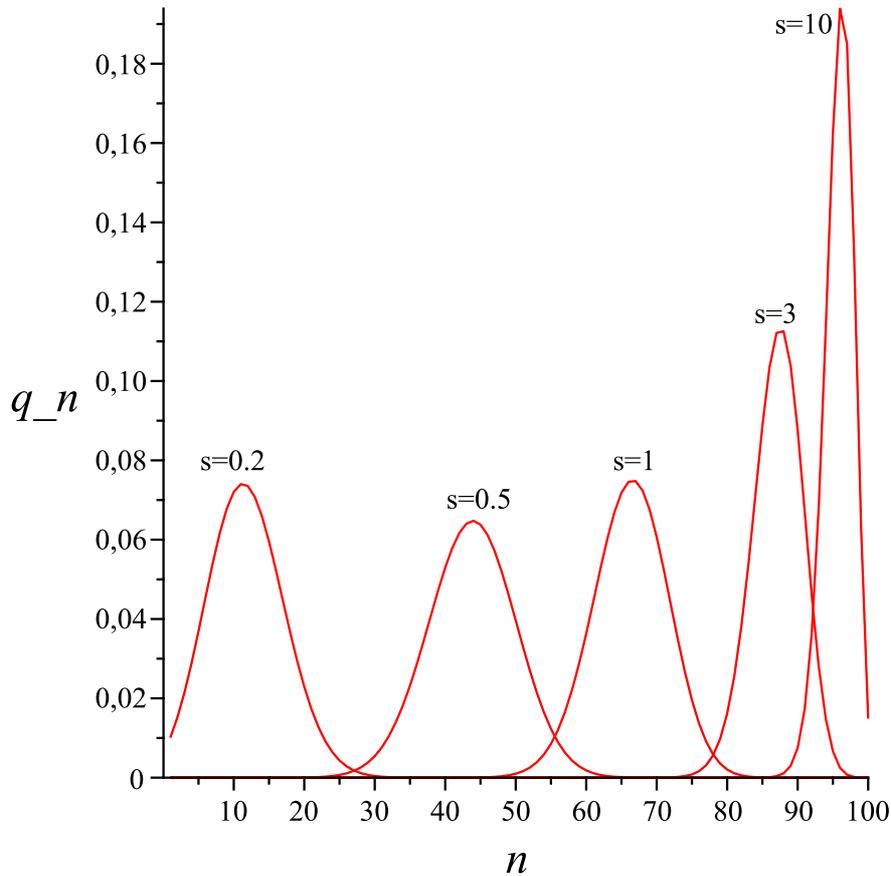,width=12cm}
\caption{The QSD for the stochastic power law logistic model with $N=100$, $R_0=5$, $\alpha=1$, and 
the 5 $s$-values 0.2, 0.5, 1, 3, and 10.}
\end{figure}

Our third observation concerns the third cumulant $\kappa_3$. 
It is useful for determining the skewness of the QSD. 
Actually, the QSD has negative skewness if $\kappa_3$ is negative, and positive skewness if 
$\kappa_3$ is positive. 
It follows from the expression \eqref{7.12} that $z_1$ is positive if $s<s_3$ and negative if 
$s>s_3$, where 
\begin{equation}
     s_3=(\alpha+1)R_0/(R_0^2+\alpha).  
\end{equation}
It is easy too see that the same intervals of $s$ lead to positive respectively negative 
skewness of the QSD if $N$ is sufficiently large.

Our brief study of the $s$-dependence of the first 3 cumulants shows that one can identify 2 cases with 
rather different behaviors. 
The first case occurs when $s<\min(s_2,s_3)$, and the second one when $s>\max(s_2,s_3)$. 
In the first case we conclude that both $\kappa_1$ and $\kappa_2$ are increasing functions of $s$, and that the QSD has positive skewness. 
In the second case we conclude as in the first case that $\kappa_1$ is increasing in $s$, 
while $\kappa_2$ now is a decreasing function of $s$ and the skewness of the QSD has changed 
from positive to negative.

We illustrate the $s$-dependence of the first three cumulants by plotting the QSDs for 5 different 
$s$-values with the constant parameter values $N=100$, $R_0=5$, and $\alpha=1$. 
In these cases we find $s_2=0.4055$ and $s_3=0.3846$.                       
Numerically determined QSDs are shown in Figure 1 for the 5 $s$-values 0.2, 0.5, 1, 3, and 10. 
The figure shows that the first cumulant (the expectation) increases as a function of $s$, 
that the second cumulant (the variance) decreases both when $s$ approaches small positive values, and when $s$ grows towards large values, and finally that the third cumulant (a measure of skewness) 
is positive for small $s$-values, and negative for large s-values. 
(For the interpretations of the plots in Figure 1, it is useful to note that the probabilities $q_n$ are positive 
for all $n$ in the state space $\{1,2,\dots,N\}$. 
Recall that the QSDs are discrete, and that the individual probabilities $q_n$ are determined 
from the plotted curves by reading off the values at each abscissa $n$.) 

The asymptotic approximations that we have given of the first three cumulants for large values of $N$     are valid for sufficiently large $N$-values. 
It is seen from the plot in Figure 1 that the requirement that $q_1$ is exponentially small is not satisfied 
for $s=0.2$ when $R_0=5$, $N=100$ and $\alpha=1$. 
We leave it as an open problem to determine how large $N$ must be as a function of $R_0$, $\alpha$ 
and $s$ to assure that $q_1$ is exponentially small in $N$.

\section{Discussion of Other Approaches}

In this section we discuss five published results that all deal with the problem in this 
paper, namely to determine approximations of the first three cumulants of the QSD 
of the stochastic power law logistic model. 
All these results are formally different from each other, and also different from 
the results that we have derived and presented in Section 7 of this paper. 
We use the powerful method of asymptotic approximations to analyze these results. 
In this discussion we are mainly limited to the case $s=1$.

The first published results are those due to Bartlett (1957) and Bartlett, Gower, ansd Leslie (1960). 
The latter results are referred to as the BGL-approximations. 
They give approximations of the first three cumulants of the QSD of the stochastic logistic 
model with $s=1$. 
The model they analyze is the first one formulated in Section 2. 
Their derivation is based on a study of the change of the stochastic variable $X(t)$, and of its seond 
and third powers, during an infinitesimal time interval.  
The resulting approximations of the first three cumulants of the QSD, which are equal to 
the mean $m=\kappa_1$, the variance $\sigma^2=\kappa_2$, and the third central moment 
$\mu_3=\kappa_3$, are given as follows in Bartlett {\it et al} (1960):
\begin{align} \label{8.1} 
   \mu         & \approx \frac{a_1-a_2}{b_1+b_2}, \\ \label{8.2}
   \sigma^2 & \approx \frac{a_1-b_1\mu}{b_1+b_2}, \\ \label{8.3}
   \mu_3      & \approx \frac{b_2-b_1}{b_1+b_2} \sigma^2.
\end{align}
In the quoted paper, these three approximations are denoted by the sign $\sim$ instead of $\approx$. 
This would indicate that the approximations have the desirable property of being asymptotic under the condition that some identified variable becomes large. 
However, no such variable appears in the derivations. 
Because of this, the derived approximations can not be claimed to be asymptotic. 
It is, however, easy to show that the BGL-approximations are closely related to asymptotic approximations of the first three cumulants. 
All one has to do is to use the reparametrization in \eqref{2.6}--\eqref{2.9}. 
This leads to the following results:
\begin{align} \label{8.4} 
   \kappa_1 & \approx x_1 N + x_2, \\ \label{8.5}
   \kappa_2 & \approx y_1 N, \\ \label{8.6} 
   \kappa_3 & \approx z_1 N, 
\end{align}
where $x_1$, $x_2$, $y_1$, $z_1$ are given in \eqref{7.1}, \eqref{7.2}, \eqref{7.4}, \eqref{7.6}.  
This shows that the approximation \eqref{8.1} of the first cumulant of the QSD is equal to the sum 
of the first two terms of its asymptotic approximation for large $N$-values, while the approximations 
\eqref{8.2} and \eqref{8.3} of the second and third cumulants of he QSD are equal to the first 
terms of the corresponding asymptotic approximations. 
All three approximations are therefore asymptotic for large $N$. 
We note, however, that the BGL-method does not by itself allow this important conclusion. 
Additional information about properties of the method will be uncovered in our study below of published 
results number three and four, referred to as BR1 and BR2. 

The second published result that we comment on was presented by Matis and Kiffe (1996) 
in case $s=1$, and by Matis, Kiffe, and Parthasarathy (1998) when $s$ is a positive integer. 
A valuable contribution of these papers is that they show that a system of ordinary differential equations (ODEs) for the first few cumulants of the unconditioned random variable $X(t)$ can be derived from 
the PDE for the cumulant generating function $K(\theta,t)$ of $X(t)$.  
As already mentioned, we use a slight variation of this approach for deriving ODEs of the first few cumulants of the conditioned random variable $X^Q(t)$. 
It turns out that the first system of ODEs (for the cumulants of $X(t)$) serves as an approximation of the second system of ODEs (for the cumulants of $X^Q(t)$) under the condition that $R_0>1$. 
The critical points of the two systems of ODEs correspond to stationary distributions of the two random 
variables $X(t)$ and $X^Q(t)$, respectively. 
The two systems of ODEs for the cumulants are not closed, in the sense that the number of 
cumulants is larger than the number of equations. 
This may appear as an undesirable property of the problem. 
Closure is clearly necessary if one wants or needs exact solutions. 
However, for our purposes it is important to realize that closure is not needed to find approximate 
solutions. 
To achieve closure, Matis {\it et al.} (1996), (1998) assume that all cumulants of sufficiently large 
order are equal to zero. 
This assumption is clearly at odds with our finding that all cumulants of the QSD are of order $O(N)$. 
It turns out, however, that it leads to good numerical approximations of the cumulants.  
This behavior can be understood from our results. 
We note for example from equations \eqref{4.15}-\eqref{4.17} that the right-hand sides of 
the ODEs for the first three cumulants when $s=4$ depend upon all seven cumulants 
$\kappa_1-\kappa_7$. 
But we note also in this case from \eqref{6.52}-\eqref{6.57} that our assumption that the cumulants 
$\kappa_4-\kappa_7$ are all $O(N)$ leads us to conclude that they have no influence 
on the first three terms of the asymptotic appoximation of $\kappa_1$, nor on the first two terms 
of the asymptotic approximation of $\kappa_2$, or on the first term of the asymptotic approximation 
of $\kappa_3$. 
The assumption that $\kappa_k=O(N)$ made by Matis and Kiffe is clearly fulfilled if $\kappa_k=0$. 
This argument shows that cumulant approximations based on cumulant closure can
lead to asymptotic approximations if they are acceptable at all, as shown in several cases by 
Matis {\it et al.} (1998). 
But the assumption that some cumulant $\kappa_k$ equals zero is incorrrect and can lead to 
errors if one studies e.g. the QSD instead of its cumulants.

The remaining three published results that we discuss here are all based on the BGL-method. 
Renshaw (2011) has used the BGL-method in an effort to derive additional terms in the approximations 
of the first three cumulants of the QSD. 
He proceeds to formulate two approaches that we refer to as BR1 and BR2, respectively, where 
the letters B and R are used to refer to Bartlett and Renshaw.
In each of these two approaches he proceeds to derive three relations between cumulants that turn out 
to be similar to relations between critical points of the system of ODEs for the first three cumulants. 
To describe his findings, we quote first from Renshaw (2011) the ODEs for the first three cumulants 
in the case $s=1$. 
They are derived from the PDE for the cumulant generating function, and do not involve 
any approximations. 
For convenience in notation we denote the time derivatives of the first three cumulants by 
$A$, $B$, and $C$, respectively. 
We get  
\begin{align} \label{8.7}
   & \kappa_1' = A =a \kappa_1 - b (\kappa_1^2+\kappa_2), \\ \label{8.8}
   & \kappa_2' = B = 2 a \kappa_2 - b(4\kappa_1 \kappa_2 + 2\kappa_3) 
       + c \kappa_1 - d(\kappa_1^2 +\kappa_2), \\ \label{8.9} 
   & \kappa_3' = C = a (\kappa_1+3\kappa_3) - b(\kappa_1^2 + 6 \kappa_1\kappa_3 + 6 \kappa_2^2 
      + \kappa_2 + 3 \kappa_4) + 3c \kappa_2 - d(6 \kappa_1 \kappa_2 + 3 \kappa_3).
\end{align}
The same equations are also found in Matis and Kiffe (1996), and, after using the reparametrization in 
\eqref{2.10}-\eqref{2.13}, in N{\aa}sell (2017).

The two approaches disussed by Renshaw (2011) are the third and fourth of the published results 
that we discuss. 
It turns out that the first two of the cumulant relations derived by Renshaw are in each of his two approaches equal to the expressions found by setting the cumulant derivatives $\kappa_1' = A$ and 
$\kappa_2' = B$ quoted above equal to zero. 
However, the third cumulant relation that he derives is different from what is found by setting 
$\kappa_3' = C$ equal to zero in each of the 2 approaches. 
In the approach BR1 the third cumulant relation is given by relation (3.5.21) in Renshaw's book, while it is given by the expression following (3.5.38) in the same book in the approach BR2. 
We use the notations $C^{(a)}=0$ and $C^{(b)}=0$ to refer to these relations. 
For the first one of these we find that $C^{(a)}$ is written as follows after 
using the reparametrization in \eqref{2.10}-\eqref{2.13}:

\begin{multline} \label{8.10}
   C^{(a)} = \mu (R_0-1) (\kappa_1^3 + 3 \kappa_1 \kappa_2 + \kappa_3) 
                   + \mu (R_0+1) \kappa_1^2 -\mu \frac{R_0-\alpha}{N} \kappa_1^3  \\ 
                - \mu \frac{R_0+\alpha}{N} (\kappa_1^4 + 6 \kappa_1^2 \kappa_2 
                + 4 \kappa_1 \kappa_3 + 3 \kappa_2^2 + \kappa_4).
\end{multline}

We proceed to derive asymptotic approximations of the first three cumulants of 
the QSD for this case. 
As in Section 6, our basic assumption is that the first four cumulants have the following 
asymptotic behaviors:
\begin{align} \label{8.11}
   \kappa_1 & = x_1 N + x_2 + \frac{x_3}{N} + \text{O}\left(\frac{1}{N^2}\right),
 \\ \label{8.12} 
   \kappa_2 & = y_1 N + y_2 + \text{O} \left(\frac{1}{N} \right), \\ \label{8.13}
   \kappa_3 & = z_1 N + \text{O}(1), \\ \label{8.14}
   \kappa_4 & = u_1 N + \text{O}(1).
\end{align}

By inserting these asymptotic approximations of the first four cumulants into 
the expressions for $A, B, C^{(a)}$ in \eqref{8.7}, \eqref{8.8},  \eqref{8.10}  
we find that asymptotic approximations for them can be written as follows:
\begin{align} \label{8.15}
    A & = A_1 N + A_2 + \frac{A_3}{N} + \text{O}\left( \frac{1}{N^2} \right), 
        \\ \label{8.16}
    B & = B_1 N + B_2 + \text{O} \left( \frac{1}{N} \right), \\ \label{8.17}
    C^{(a)} & = C^{(a)}_3 N^3 + C^{(a)}_2 N^2 + C^{(a)}_1 N + \text{O}(1),
\end{align}
where

\begin{align} \label{8.18}
   A_1 & = \mu (R_0-1) x_1 -\mu (R_0+\alpha) x_1^2, \\ \label{8.19}
   A_2 & = \mu (R_0-1) x_2 - \mu (R_0+\alpha) (2 x_1 x_2 + y_1), \\ \label{8.20}
   A_3 & = \mu (R_0-1) x_3 - \mu (R_0+\alpha) (2 x_1 x_3 + x_2^2 +y_2), \\ \label{8.21} 
   B_1 & = 2\mu(R_0-1) y_1 + \mu(R_0+1) x_1 -\mu (R_0-\alpha) x_1^2 
               - 4\mu (R_0+\alpha) x_1 y_1 , \\ \label{8.22} 
   B_2 & = 2\mu (R_0-1) y_2  + \mu (R_0+1) x_2 - \mu (R_0-\alpha) (2 x_1 x_2 + y_1) \\ \notag
& \phantom{he}    - \mu (R_0+\alpha) (4 x_1 y_2 + 4 x_2 y_1+ 2 z_1), \\ \label{8.23}
   C^{(a)}_3 & =\mu (R_0-1) x_1^3 -\mu (R_0+\alpha) x_1^4, \\ \label{8.24}
   C^{(a)}_2 & = \mu (R_0-1) (3 x_1^2 x_2 + 3 x_1 y_1) + \mu (R_0+1) x_1^2 
- \mu (R_0-\alpha) x_1^3 \\ \notag
 & \phantom{he} - \mu (R_0+\alpha) (4 x_1^3 x_2 + 6 x_1^2 y_1),   \\ \label{8.25}
   C^{(a)}_1 & = \mu (R_0-1) (3 x_1^2 x_3 + 3 x_1 x_2^2 + 3 x_1 y_2 + 3 x_2 y_1 + z_1) 
\\ \notag
 & \phantom{he} + 2 \mu (R_0+1) x_1 x_2 - 3 \mu (R_0-\alpha) x_1^2 x_2  \\ \notag
 & \phantom{he} - \mu (R_0 + \alpha) (4 x_1^3 x_3 + 6 x_1^2 x_2^2 + 6 x_1^2 y_2 
          + 12 x_1 x_2 y_1 + 4 x_1 z_1 + 3 y_1^2). 
\end{align}

The basic mathematical problem at this point is to determine the six coefficients 
$x_1, x_2, x_3, y_1, y_2, z_1$ so that the following three conditions are satisfied: 
\begin{align} \label{8.26}
   & A = \text{O} \left( \frac{1}{N^2}  \right), \\  \label{8.27}
   & B = \text{O} \left( \frac{1}{N}  \right), \\  \label{8.28}
   & C^{(a{)}} = \text{O} \left( 1 \right).
\end{align}

These conditions are satisfied by setting the eight expressions \eqref{8.18} - \eqref{8.25} equal to zero.  
We proceed to solve the resulting eight equations for the six unknown coefficients 
$x_1$, $x_2$, $x_3$, $y_1$, $y_2$, $z_1$. 
At this point it appears that there are more equations than unknowns, but the apparent 
problem that this causes will readily be solved. 
The equation $A_1=0$ is first solved for $x_1$, and the spurious solution 
$x_1=0$ is excluded. 
As above, as soon as the solution is found for any of the coefficients $x_1$, $x_2$, $x_3$, 
$y_1$, $y_2$, $z_1$, it is inserted into the remaining unsolved equations. 
The equation $B_1=0$ is then solved for $y_1$, and the equation $A_2=0$ is 
solved for $x_2$. 
After this we find that the two expressions $C^{(a)}_3$ and $C^{(a)}_2$ can both 
be determined from these values, and that they both are equal to zero. 
Thus, the number of equations is reduced to be equal to the number of unknown 
coefficients.   
In the last steps we use the equation $A_3=0$ to express $x_3$ as a function 
of $y_2$, the equation $B_2=0$ to solve for $z_1$ as a function of $y_2$, and 
finally the equation $C_1^{(a)}=0$  to solve for $y_2$.   
The result of these evaluations is that the six coefficients $x_1$, $x_2$, $x_3$, $y_1$, 
$y_2$, $z_1$ can be expressed as follows as functions of the model parameters 
$R_0$ and $\alpha$:
\begin{align} \label{8.29}
   x_1 & = \frac{R_0-1}{R_0+\alpha}, \\ \label{8.30}
   x_2 & = - \frac{(\alpha+1)R_0}{(R_0+\alpha)(R_0-1)}, \\ \label{8.31}
   x_3 & = \frac{(\alpha+1) R_0 (R_0^2+\alpha - 5 (\alpha+1) R_0)}
              {(R_0+\alpha) (R_0-1)^3},             \\ \label{8.32}
   y_1 & = \frac{(\alpha+1) R_0}{(R_0+\alpha)^2}, \\ \label{8.33}
   y_2 & = - \frac{(\alpha+1) R_0 (R_0^2 + \alpha - 4 (\alpha+1) R_0) }
              {(R_0+\alpha)^2 (R_0-1)^2 }, \\ \label{8.34}
   z_1 & = \frac{(\alpha+1) R_0 (R_0^2+\alpha - 3 (\alpha+1) R_0) }
             {(R_0+\alpha)^3 (R_0-1) }.
\end{align}

Comparisons with \eqref{7.1}-\eqref{7.6} show that the expressions for the 3 coefficients 
$x_1$, $x_2$, and $y_1$ are equal to the corresponding expressions derived by N{\aa}sell (2017) 
using our preferred method where asymptotic approximations of the cumulants are derived from 
the ODEs for the cumulants without attempting to close the equations.  
The comparisons also show that the expressions for the remaining three coefficients $x_3$, $y_2$, 
and $z_1$ disagree between our approach and BR1.  
This means that the two approaches lead to different approximations of all three cumulants. 
Since both of them cannot be right, and since we claim that our derivations are correct, 
we conclude that the BR1 results are incorrect. 
The reason for the discrepancy is that one or more of the various approximation steps that are taken in using the BGL-method brings in errors. 
Renshaw (2011) appears to share this view in his brief discussion of this issue after his formula (3.5.21). 
We conclude that the BGL-method brings in errors of unknown magnitude, and therefore is unsuitable for deriving approximations of the cumulants of the QSD beyond those derived by Bartlett {\it et al.} (1960).  
Support for this claim is given by the numerical evaluations reported in Table 5. 
They indicate that the error terms of the approximations of the cumulants  
$\kappa_1, \kappa_2, \kappa_3$, derived in N{\aa}sell (2017), are of the orders O(1/$N^2$),  
O(1/$N$), O(1), respectively, just as expected. 
But we find also that the error terms of the approximations derived via the BR1 approach are 
larger than this by an order of magnitude. 
The BR1 approximation of $\kappa_3$ in this case is useless, since the error terms 
actually grow with $N$.

\begin{table}[h]
  \begin{center}
    \begin{tabular}{ | c | c | c | c | c | }      \hline
  Approx  & Cumulant & N=100                      & N=200                & N=400             \\ \hline                  
     N      & $\kappa_1$   & $-697*10^{-7}$   & $-171*10^{-7}$  & $-42*10^{-7}$ \\
     N      &  $\kappa_2$  & $445*10^{-5}$    & $219*10^{-5}$   & $108*10^{-5}$ \\
     N      & $\kappa_3$   & $-226*10^{-3}$   & $-222*10^{-3}$  & $-220*10^{-3}$ \\  \hline
    BR1   & $\kappa_1$   & $-311*10^{-5}$   & $-154*10^{-5}$  & $-76*10^{-5}$ \\
    BR1   & $\kappa_2$   & $253*10^{-3} $   & $251*10^{-3} $  & $250*10^{-3}$ \\
    BR1   & $\kappa_3$   & -21                      & -41                     & -82                     \\ \hline
     BB    & $\kappa_1$   & $-359*10^{-5}$   &  $-178*10^{-5}$ & $-88*10^{-5}$ \\                  
     BB    & $\kappa_2$   &   $292*10^{-3}$  & $290*10^{-3}$   & $289*10^{-3}$ \\
     BB    & $\kappa_3$   & $-10.2$               & $-20.3$               & $-40.3$             \\ \hline
      \end{tabular}
    \vskip 4mm
    \caption{Numerical evaluations of the error terms of 3 different 
 approximations of the first 3 cumulants of the QSD of the stochastic Verhulst 
 logistic model. 
 Results are shown for $R_0=10$, $\alpha=1$, $s=1$,  and the $N$-values 100, 
 200, and 400. 
 Approximation N uses the results presented in this paper, while Approximation  BR1 uses the 
 first version of the Bartlett-Renshaw approach, and Approximation BB is taken from Section 4 
 of the paper by Bhowmick {\it et al.}
 The table indicates that the error terms of $\kappa_1, \kappa_2, \kappa_3$ in  Approximation N 
 are reduced by approximately 4, 2, and 1, respectively, for  each doubling of $N$, while in 
 approximations BR1 and BB the error terms of  $\kappa_1$ and $\kappa_2$ are reduced by 
 approximately 2 and 1, respectively,  for each doubling of $N$, and the error terms of 
 $\kappa_3$ are seen to grow  with $N$. 
 The derivations of the numerical results here are documented in N{\aa}sell (2018).}
  \end{center}
\end{table}

We turn now to a consideration of the second version BR2 of the Bartlett-Renshaw approach. 
Renshaw (2011) uses it to derive a new expression for the third cumulant relation. 
It is given after relation (3.5.38) in his book. 
After writing $\mu_4 = \kappa_4 + 3 \kappa_2^2$ and using the reparametrization \eqref{2.10}-\eqref{2.13} it can be written $C^{(b)}=0$, where $C^{(b)}$ is as follows:
\begin{multline} \label{8.35} 
   C^{(b)} = \mu (R_0-1) (\kappa_1 \kappa_2 + \kappa_3) + \mu (R_0+1) \kappa_2
                   - \mu \frac{R_0-\alpha}{N} (2 \kappa_1 \kappa_2 + \kappa_3)  \\
                - \mu \frac{R_0+\alpha}{N} (\kappa_1^2 \kappa_2 + 2 \kappa_1 \kappa_3
                    + 3 \kappa_2^2 + \kappa_4).
\end{multline}

We proceed to derive asymptotic approximations of the first three cumulants of the 
QSD for this case, using the three cumulant relations $A=0$, $B=0$, and 
$C^{(b)}=0$.  
As above, we assume that the first four cumulants have the following 
asymptotic behaviors:
\begin{align} \label{8.36} 
    \kappa_1 & = x_1 N + x_2 + \frac{x_3}{N} +  \text{O}\left(\frac{1}{N^2} 
 \right), \\ \label{8.37}
    \kappa_2 & = y_1 N + y_2 + \text{O}\left( \frac{1}{N} \right), \\ 
 \label{8.38}
    \kappa_3 & = z_1 N + \text{O}(1), \\  \label{8.39}
    \kappa_4 & = u_1 N + \text{O}(1). 
\end{align}

By inserting these asymptotic approximations of the first four cumlants into 
the expressions for $A$, $B$, and $C^{(b)}$, we find that asymptotic 
approximations for them can be written as follows:
\begin{align} \label{8.40}
&   A = A_1 N + A_2 + \frac{A_3}{N} + \text{O}\left( \frac{1}{N^2} \right),  
\\ \label{8.41} 
&   B = B_1 N + B_2 + \text{O}\left( \frac{1}{N} \right),  \\ \label{8.42}
&   C^{(b)} = C^{(b)}_2 N^2 + C^{(b)}_1 N + \text{O}(1).
\end{align}
Expressions for the five quantities $A_1, A_2, A_3, B_1, B_2$ are given in 
\eqref{8.18}-\eqref{8.22}, 
while $C^{(b)}_2$ and $C^{(b)}_1$ are equal to  
\begin{align} \label{8.43} 
&   C^{(b)}_2 = \mu (R_0-1) x_1 y_1 - \mu (R_0+\alpha) x_1^2 y_1,  \\  \label{8.44} 
&   C^{(b)}_1 = \mu (R_0-1) (x_1 y_2 + x_2 y_1 + z_1) + \mu (R_0+1) y_1 
                         - 2 \mu (R_0 - \alpha) x_1 y_1 \\ \notag
&   \phantom{hejhej} - \mu (R_0+\alpha) 
   \left(x_1^2 y_2 + 2 x_1 x_2 y_1 + 2 x_1 z_1 + 3 y_1^2\right).
\end{align}    

We set these 7 expressions equal to zero, and solve the resulting seven equations 
for the six unknown coefficients $x_1, x_2, x_3, y_1, y_2, z_1$. 
The value of $x_1$ is found by solving the equation $A_1=0$ and excluding the 
spurious solution $x_1=0$. 
The values of $y_1$ and $x_2$ are then found by first solving the equation 
$B_1=0$ for $y_1$, and then solving the equation $A_2=0$ for $x_2$. 
The values of $x_1$, $x_2$, and $y_1$ are the same as the ones found in 
\eqref{8.29}, \eqref{8.30}, \eqref{8.32}. 
Furthermore, we insert these values of $x_1, x_2, y_1$ into \eqref{8.43} 
and \eqref{8.44}.
This shows that $C^{(b)}_2=0$ and that therefore the number of equations available for 
solving for the six coefficients is reduced from seven to six. 
We also get 
\begin{equation} \label{8.45} 
    C^{(b)}_1 = -(R_0-1) z_1 
                   - \frac{(\alpha+1) R_0 (R_0-\alpha) (R_0-1)}{(R_0+\alpha)^3}. 
\end{equation}
To determine the three coefficients $x_3$, $y_2$, $z_1$ we solve the equation 
$C^{(b)}_1=0$ for $z_1$, the equation $B_2=0$ for $y_2$, and the equation 
$A_3=0$ for $x_3$. 
The rather surprising result is that the values of the six coefficients $x_1$, $x_2$, 
$x_3$, $y_1$, $y_2$, $z_1$ are found to be the same as in \eqref{7.1}-\eqref{7.6}. 
This means that even though the explicit results from the BR2 method are different from the 
results from our preferred method, the asymptotic approximations agree.

In his development of BR2, Renshaw (2011) uses the equation $C^{(b)}_1=0$ to derive the relation
\begin{equation} \label{8.46} 
     \kappa_3 \approx -\frac{R_0-\alpha}{R_0+\alpha} \kappa_2. 
\end{equation} 
By using the reparametrization in \eqref{2.11} and \eqref{2.13} we find that this relation 
was shown to hold already by Bartlett {\it et al.} (1960). 
It is seen from \eqref{7.4} and \eqref{7.6} that this relation is asymptotic for large $N$. 
However, the arguments used by Renshaw in his derivation can be criticized. 
As also pointed out by Bhowmick {\it et al.} (2016), it does not make sense to 
assume that $\kappa_1$ is large in comparison with $a_1, b_1, a_2, b_2$, since 
these four parameters are rates whose values depend on the unit of time, 
which is arbitrary, while $\kappa_1$ is independent of the time unit.
It is also unrealistic to assume that $\kappa_2$ is small compared with 
$\kappa_1$, since it contradicts our finding that all cumulants are $O(N)$.

Renshaw's two efforts to extend the BGL-approximation beyond what was derived by 
Bartlett {\it et al.} (1960) leads to two different conclusions. 
In one case (BR2) the results are correct asymptotically, while they are not in the other case (BR1). 
This is enough to cnclude that the BGL-method cannot be trusted to give correct results 
wthout further investigations.

We turn now to consider the results reported by Bhowmick {\it et al.} (2016). 
They work with the parameter space that is associated with the second of the model 
formulations in Section 2. 
Actually, they study a more general model than the one that we are concerned 
with here. 
The population birth rate $\lambda_n$ in their model is given by 
\begin{equation} \label{8.47} 
   \lambda_n =  \mu R_0 \left(1 - \alpha_1 \left(\frac{n}{N}\right)^\beta \right) n^\delta,                                  
\end{equation}
while their population death rate $\mu_n$ equals
\begin{equation} \label{8.48}
     \mu_n = \mu \left(1 + \alpha_2 \left(\frac{n}{N}\right)^\beta \right) n^\delta.
\end{equation} 
They do not require $\beta = s$ to be an integer. 
We study their model in the special case where $\beta=\delta=1$. 
To agree with our model formulation we furthermore put $\alpha_1=1$ and 
$\alpha_2=\alpha$.

By using results derived by Bhowmick {\it et al.} (2016) and reported in 
Section 4 of their paper, we find that the mean $\kappa_1$ equals 
\begin{equation} \label{8.49}
    \kappa_1 = \frac{R_0-1}{R_0+\alpha} N \frac{1}{1 +H/N},
\end{equation}
where 
\begin{equation} \label{8.50}
    H = \frac{(\alpha+1) R_0}{(R_0-1)^2}.
\end{equation}
By including 3 terms in the asymptotic approximation of $\kappa_1$ in 
\eqref{8.49}, we find  that 
\begin{equation} \label{8.51} 
    \kappa_1 = \frac{R_0-1}{R_0+\alpha} N 
                      - \frac{(\alpha+1)R_0}{(R_0+\alpha)(R_0-1) }
                  + \frac{(\alpha+1)^2 R_0^2}{(R_0+\alpha)(R_0-1)^3} \frac{1}{N} 
                        + \text{O} \left( \frac{1}{N^2}  \right).
\end{equation}

From relation (25) in Bhowmick {\it et al.} we find that the variance 
$\kappa_2$ equals
\begin{equation} \label{8.52}
     \kappa_2 = \kappa_1^2 \frac{H}{N}.
\end{equation} 
By including two terms in the asymptotic approximation of $\kappa_2$, we find 
that
\begin{equation} \label{8.53} 
   \kappa_2 = \frac{(\alpha+1) R_0}{(R_0+\alpha)^2}N 
                    - \frac{2(\alpha+1)^2 R_0^2}{(R_0+\alpha)^2 (R_0-1)^2}
                    + \text{O}\left(\frac{1}{N} \right).
\end{equation}

Bhowmick {\it et al.} use four quantities $A$, $B$, $C$, and $D$ in a formula 
that determines the third cumulant $\kappa_3$ as follows:
\begin{equation} \label{8.54}
    \kappa_3 = - \frac{A+B}{C+D} \kappa_2.
\end{equation}
These four quantities are defined in terms of the parameters that are used to describe the 
first of the two model formulations in Section 2. 
After reparametrization they can be expressed as follows: 
\begin{align} \label{8.55} 
&   A = a - b \kappa_1 = \mu(R_0-1) - \mu\frac{R_0+\alpha}{N}\kappa_1, \\  \label{8.56}
&   B = \frac{c}{\kappa_1} - 2d = \mu\frac{R_0+1}{\kappa_1} - 2\mu \frac{R_0-\alpha}{N}, \\  \label{8.57} 
&   C = \frac{a}{\kappa_1} - 2b = \mu\frac{R_0-1}{\kappa_1} - 2\mu \frac{R_0+\alpha}{N},  \\  \label{8.58}
&   D = - \frac{d}{\kappa_1} = - \mu \frac{R_0-\alpha}{N\kappa_1}.
\end{align}

Thus, all four of these quantities are determined as functions of $\kappa_1$. 
We use the asymptotic approximation of $\kappa_1$ in \eqref{8.51} to determine 
one-term asymptotic approximations of each of $A, B, C, D$. 
Two terms of the asymptotic approximation of $\kappa_1$ is needed for $A$, while one term 
suffices for the remaining 3 quantities. 
The results are 
\begin{align} \label{8.59}
&   A = \mu \frac{(\alpha+1)R_0}{R_0-1} \frac{1}{N}+ \text{O}\left(\frac{1}{N^2} 
\right), \\ \label{8.60}
&   B = -\mu \frac{R_0^2+\alpha - 3(\alpha+1)R_0}{R_0-1} \frac{1}{N}
+\text{O} \left(\frac{1}{N^2} \right), \\ \label{8.61}
&   C = - \mu (R_0+\alpha)\frac{1}{N} +\text{O} \left(\frac{1}{N^2} \right), 
\\ \label{8.62}
&   D = -\mu \frac{R_0^2-\alpha^2}{R_0-1} \frac{1}{N^2} + \text{O} \left(\frac{1}{N^3} \right).
\end{align}    
It follows from this that the one-term asymptotic approximation of 
$(A+B)/(C+D)$ equals
\begin{equation} \label{8.63}
   \frac{A+B}{C+D} = \frac{R_0^2+\alpha - 4(\alpha+1) R_0}{(R_0+\alpha) (R_0-1)} 
                                 + \text{O}\left(\frac{1}{N} \right). 
\end{equation}
By inserting this result and the one-term asymptotic approximation of 
$\kappa_2$ from \eqref{8.53} into \eqref{8.54}, we get
\begin{equation} \label{8.64}
    \kappa_3 = - \frac{(\alpha+1) R_0 (R_0^2 + \alpha - 4 (\alpha+1)R_0)}    
{(R_0+\alpha)^3 (R_0-1)} N + \text{O}(1).
\end{equation}

The resulting asymptotic approximations of the first 3 cumulants are found in 
\eqref{8.51}, \eqref{8.53}, and \eqref{8.64}. 
Comparisons with the results derived by using our preferred method and found from 
\eqref{7.1}-\eqref{7.6} show that the two approaches lead to different results. 
We claim that our results are correct, and that the results given by Bhowmick 
{\it et al.} are incorrect.  
An independent suppport for this conclusion is found from the numerical evaluations 
of the error terms of the results presented by Bhowmick {\it et al.} and given in 
Table 5. 
We see here that the error term for the BB approximation of $\kappa_1$ is reduced 
by the factor 2 for each doubling of $N$. 
This is typical of an error that is of magnitude of $\text{O}(1/N)$. 
In addition, the error term for the BB approximation of $\kappa_2$ is found to 
be approximately independent of $N$, which indicates that it is of the magnitude 
of $\text{O}(1)$. 
The error term for the BB approximation of $\kappa_3$ is, finally, found to be 
approximately proportional to $N$, which makes the corresponding 
approximation useless.  
The reason for the incorrect results arrived at by Bhowmick {\it et al.} is  that one 
or more of the approximation steps taken in the application of BGL method brings in errors. 
This is similar to the reason for incorrectness of the BR1 result, commented on above.

\section{Concluding Comments}

The method that we have used here to derive asymptotic approximations of the 
first few cumulants of the QSD of the stochastic power law logistic model 
emphasizes the importance of the second model formulation in Section 2. 
It gives access to two parameters that are basic for our study, namely the 
maximum population size $N$ and the threshold parameter $R_0$. 
We have shown that the condition $R_0>1$ for large $N$ is required for the 
results that we have presented here, while both the moment closure method 
and the 3 approaches BR1, BR2, and BB based on the BGL method have 
the serious weaknesses that they can not produce conditions for validity of the 
approximations that are derived. 
In addition we note that magnitudes of approximation errors are easy to establish 
in our method, as they are for any asymptotic approximation, while they can not 
be produced in the moment closure method, nor in approaches based on the 
BGL method. 
We note furthermore that spurious solutions that could require large efforts 
to eliminate appeared in early studies based on moment closure, while the 
spurious solutions that appear in the method that we use here are easy to 
identify and eliminate.  
Our results show that the number of spurious solutions is equal to the 
parameter $s$ whenever $s$ is a positive integer. 
As a further comment we note that the dependence of the new approximations on the model 
parameters is explicit in our approach, while they are unknown in results based on 
moment closure and in BR1 and BR2.

We have found that our method based on the second model formulation of Section 2, 
and followed by a search for asymptotic approximations, provides a powerful 
approach for determination of approximations of the first few cumulants of the QSD 
for the power law logistic model. 
We have also shown that the method of determining asymptotic approximations 
can be used to study other approaches to the same problem.

We conclude that the method of determining ODEs for the first few cumulants of the QSD 
introduced by Matis and Kiffe (1996) is preferred over the BGL-method for deriving relations 
between the cumulants. 
We emphasize the obvious fact that whenever exact solutions of a mathematical problem are difficult 
to establish, then one should search for approximations. 
Furthermore, approximation methods are then definitely preferred that give both condiitions 
for validity of the approximations and magnitudes of approximation errors.  
Because of this, we conclude that our method for determining asymptotic approximations 
of the first few cumulants is preferred over any method that is based on moment closure.



\begin{thebibliography}{99}

\bibitem{B1994}
R.B. Banks,
Growth and Diffusion Phenomena, 
Springer, Berlin, Heidelberg (1994).

\bibitem{B1957}
M.S. Bartlett,
On theoretical models for competitive and predatory biological systems,
Bimetrika, {\bf 44}, 27--42, 1957. 


\bibitem{BGL1960}
M.S. Bartlett, J.C. Gower, and P.H. Leslie,
A comparison of theoretical and empirical results for some stochastic
population models,
Biometrika, {\bf 47}, 1--11 (1960).

\bibitem{BBRB2016}
A. R. Bhowmick, S. Bandyopadhyay, S. Rana, and S. Bhattacharya, 
A simple approximation of moments of the quasi-equilibrium distribution of an 
extended theta-logistic model with non-integer powers, 
Math. Biosci., {\bf 271}, 96--112 (2016). 



\bibitem{FKMP1995}
P. Ferrari, H. Kesten, S. Martínez, and P. Pico, 
Existence of quasi-stationary distributions. A renewal dynamic approach, 
Ann. Probab. {\bf 23}, 501--521 (1995).


\bibitem{K2016}
C. Kuehn,
Moment closure - A brief review, 
In Control of Self-Organizing Nonlinear Systems, Springer, 2016, 
arXiv:1505.02190


\bibitem{MK1996}
J.H. Matis and T.R. Kiffe, 
On approximating the moments of the equilibrium distribution of a stochastic 
logistic model, 
Biometrics, {\bf 52}, 980--991 (1996). 
	

\bibitem{MKP1998}
J.H. Matis, T.R. Kiffe, and P.R. Parthasarathy,
On the cumulants of population size for the stochastic power law logistic model,
Theor. Pop. Biol., {\bf 53}, 16--29 (1998).  



\bibitem{N2001a}
I. N{\aa}sell, 
Extinction and quasi-stationarity in the Verhulst logistic model, 
J. Theor. Biol., {\bf 211}, 11--27, 2001a. 

\bibitem{N2001b}
I. N{\aa}sell, 
Extinction and quasi-stationarity in the Verhulst logistic model: 
With derivations of mathematical results.
http://people.kth.se/~ingemar/forsk/verhulst/verhulst.html, (2001b).  

\bibitem{N2011}
I. N{\aa}sell, 
Extinction and Quasi-stationarity in the Stochastic Logistic SIS Model, 
Springer Lecture Notes in Mathematics, Vol {\bf 2022}, Berlin, Heidelberg (2011). 

\bibitem{N2017}
I. N{\aa}sell, 
An alternative to moment closure, 
Bull. Math. Biol., {\bf 79}, Issue 9, 2088--2108,  (2017).


\bibitem{N2018}
I. N{\aa}sell,
Maple Worksheets for the study of the cumulants of the stochastic power law 
logistic model,
https://people.kth.se/$\sim$ingemar/PLawLogistic, (2018).



\bibitem{R2011}
E. Renshaw,
Stochastic Population Processes: Analysis, Approximations, Simulations,
Oxford University Press, Oxford (2011).  


 
\bibitem{TW2002}
A. Tsoularis and J. Wallace,
Analysis of logistic growth models,
Math. Biosci., {\bf 179}, 21-55 (2002). 
 
\bibitem{V1838}
P.F. Verhulst, 
Notice sur la loi que la population suit dans son accroisement,
Corr. Math. Phys. {\bf X}, 113--121 (1838).  

\bibitem{W57}
P. Whittle, 
On the use of the normal approximation in the treatment of stochastic 
processes,
J. Roy. Statist. Soc., Ser. B 19, 268--281 (1957).  


\end{thebibliography}
\end{document}